\documentclass[11pt]{article}

\usepackage{graphicx}
\usepackage{transparent}
\usepackage{geometry}
\usepackage{amsmath}
\usepackage{amsfonts}
\usepackage{amsthm}
\usepackage{amssymb}
\usepackage{mathdots}
\usepackage{latexsym}
\usepackage{mathrsfs}
\usepackage{subfig}
\usepackage{enumerate}
\usepackage{rotating}
\usepackage{multirow}
\usepackage{hyperref}
\usepackage{pgfplots}
\usepackage{lineno, blindtext}
\usepackage{cancel}
\usepackage{wrapfig}
\usepackage{float}
\usepackage{multicol} 
\usepackage{graphicx}
\usepackage{dsfont}
\usepackage{yfonts}
\usepackage{cite}
 
\newtheorem{theorem}{Theorem}[section]
\newtheorem{lemma}[theorem]{Lemma}
\newtheorem{proposition}[theorem]{Proposition}
\newtheorem{corollary}[theorem]{Corollary}
\newtheorem{definition}[theorem]{Definition}
\newtheorem{remark}[theorem]{Remark}
\newtheorem{question}[theorem]{Question}

\newcommand{\C}{\mathcal{C}}
\newcommand{\Mod}{{\text{Mod}}} 
\newcommand{\pf}{\pitchfork}
\newcommand{\wt}{\widetilde}
\newcommand{\Hy}{\mathbb{H}^2}
\newcommand{\ol}{\overline}

\newcommand{\hgline}[2]{
\pgfmathsetmacro{\thetaone}{#1}
\pgfmathsetmacro{\thetatwo}{#2}
\pgfmathsetmacro{\theta}{(\thetaone+\thetatwo)/2}
\pgfmathsetmacro{\phi}{abs(\thetaone-\thetatwo)/2}
\pgfmathsetmacro{\close}{less(abs(\phi-90),0.0001)}
\ifdim \close pt = 1pt
    \draw[black] (\theta+180:1) -- (\theta:1);
\else
    \pgfmathsetmacro{\R}{tan(\phi)}
    \pgfmathsetmacro{\distance}{sqrt(1+\R^2)}
    \draw[black] (\theta:\distance) circle (\R);
\fi
}

\newcommand{\hglinewhite}[2]{
\pgfmathsetmacro{\thetaone}{#1}
\pgfmathsetmacro{\thetatwo}{#2}
\pgfmathsetmacro{\theta}{(\thetaone+\thetatwo)/2}
\pgfmathsetmacro{\phi}{abs(\thetaone-\thetatwo)/2}
\pgfmathsetmacro{\close}{less(abs(\phi-90),0.0001)}
\ifdim \close pt = 1pt
    \draw[black] (\theta+180:1) -- (\theta:1);
\else
    \pgfmathsetmacro{\R}{tan(\phi)}
    \pgfmathsetmacro{\distance}{sqrt(1+\R^2)}
    \fill[white] (\theta:\distance) circle (\R);
\fi
}

\begin{document}

\title{Non-existence of boundary maps for  some \\ hierarchically hyperbolic spaces}
\author{Sarah C. Mousley}
\date{}

\maketitle

\abstract{We provide negative answers to questions posed by Durham, Hagen, and Sisto on the existence of boundary maps for some hierarchically hyperbolic spaces, namely maps from right-angled Artin groups to mapping class groups. We also prove results on existence of boundary maps for free subgroups of mapping class groups.}

\section{Introduction}

Let $\Gamma$ be a finite graph with vertex set $V(\Gamma)=\{s_1, \ldots, s_k\}$. The right-angled Artin group determined by $\Gamma$, denoted by $A(\Gamma)$, is the group with the following presentation:
\[ A(\Gamma)=\langle s_1, \ldots, s_k : [s_i, s_j]=1 \Leftrightarrow s_is_j \text{ is an edge in } \Gamma \rangle.\]

Let $S=S_{g,n}$ be a connected, oriented surface of genus $g$ with $n$ punctures, and let $\Mod(S)$ denote the mapping class group of $S$.
Clay, Leininger, and Mangahas \cite{CLM} and Koberda \cite{Koberda} construct ``nice'' embeddings of right-angled Artin groups to mapping class groups. In \cite{BHS} and \cite{BHSII}, a geometric structure called a hierarchically hyperbolic space (HHS) was introduced. Important examples of spaces that are HHS's include mapping class groups of surfaces and right-angled Artin groups.
In \cite{DHS} Durham, Hagen, and Sisto construct a boundary for hierarchically hyperbolic spaces (see Section \ref{sec:background}). In that paper, the  authors ask the following question, motivated by a desire to develop a notion of geometrically finite subgroups of mapping class groups.

\begin{question} \label{mainq}
Let $A(\Gamma)$ be a right-angled Artin group embedded in $\Mod(S)$ in the sense of either Clay, Leininger, and Mangahas \cite{CLM} or Koberda \cite{Koberda}. Does the embedding $ A(\Gamma) \rightarrow \Mod(S)$ extend continuously to an injective map $\partial A(\Gamma) \rightarrow  \partial Mod(S)$?
\end{question}

We prove that in general the answer to this question is no by providing, for each type of embedding, an explicit example where the embedding does not extend continuously.

\begin{theorem} \label{theoremCLMKintro}
There exists a surface $S$, a right-angled Artin group $\Gamma$, a Clay, Leininger, and Mangahas  embedding $\phi: A(\Gamma) \rightarrow \Mod(S)$, and a Koberda embedding $\phi':A(\Gamma) \rightarrow \Mod(S)$  such that, regardless of the HHS structure on $A(\Gamma)$, neither $\phi$ nor $\phi'$ extends continuously to a map $\partial A(\Gamma) \rightarrow  \partial \Mod(S)$.
\end{theorem} 

We also prove the following result which gives a complete characterization of Koberda embeddings of free groups, which send all generators to powers of Dehn twists, that have continuous extensions. 

\begin{theorem} \label{theorem:IntroExistenceFill}  Let $\{\alpha_1, \ldots, \alpha_k\}$ be a collection of pairwise intersecting curves in $S$  and  $\Gamma$ the graph with $V(\Gamma)=\{s_1, \ldots, s_k\}$ and no edges. 
For sufficiently large $N$, the  homomorphism \[\phi: A(\Gamma) \rightarrow \Mod(S) \hspace{10pt} \text{ defined by } \hspace{10pt} \phi(s_i)= T_{\alpha_i}^{N}  \text{ for all } i \] is injective by the work of Koberda \cite{Koberda}. Moreover, $\phi$ extends continuously  to a map $ \partial A(\Gamma) \rightarrow \partial \Mod(S)$ if and only if $\{\alpha_1, \ldots ,\alpha_k\}$  pairwise fill $S$, where $A(\Gamma)$ is equipped with any HHS structure. 
\end{theorem}

In fact, we prove something stronger than Theorem \ref{theorem:IntroExistenceFill}. We prove a non-existence result (Theorem \ref{theorem:NofillNoextend})  for a class of Koberda embeddings of right-angled Artin groups
 that are not necessary free groups. We also prove an existence result (Theorem \ref{theorem:existenceFill}) for a class of embeddings of free groups
that includes the Koberda embeddings described in Theorem \ref{theorem:IntroExistenceFill} as well as a class of Clay, Leininger, Mangahas embeddings. 

In Section \ref{sec:background} we will recall relevant definitions and theorems and introduce notation. Section \ref{sec:lemma} will establish a handful of lemmas that will be used for proving Theorem \ref{theoremCLMKintro}. Section \ref{sec:CLMRAAG} is devoted to proving Theorem \ref{theoremCLMKintro} for a Clay, Leininger, Mangahas embedding, and in Section \ref{sec:KRAAG} we prove Theorem \ref{theoremCLMKintro} for a Koberda embedding. Using similar techniques, we then  prove that a more general class  of Koberda embeddings of right-angled Artin groups do not extend continuously (Theorem \ref{theorem:NofillNoextend}), which will imply one direction of Theorem \ref{theorem:IntroExistenceFill}. In Section \ref{section:existence} we will prove Theorem \ref{theorem:existenceFill}, which will imply the other direction of Theorem \ref{theorem:IntroExistenceFill}.  

\

\noindent \textbf{Remark:} Koberda \cite{Koberda} proved that both types of embeddings we discuss are injective. We call the embeddings that send generators of our right-angled Artin group to mapping classes that are pseudo-Anosov on subsurfaces Clay, Leininger, Mangahas embeddings primarily to distinguish the two types, but also to emphasize that these types of embeddings have nice geometric properties (see Theorem \ref{CLMemb}). 

\

{\bf Acknowledgments:} The author was supported by the Department of Defense (DoD) through the National Defense Science and Engineering Graduate Fellowship (NDSEG) Program and by a Research Assistantship through NSF Grant number DMS-1510034. The author would like to thank her PhD advisor Chris Leininger for   his numerous insights which inspired a great deal of this paper and also for his patience and consistent support. The author would also like to thank Mark Hagen and Matt Durham for helpful conversations. 

\section{Background} \label{sec:background}
In this section, we recall some needed definitions and theorems.

\

\noindent {\bf Notation:} Let $f, g: X \rightarrow \mathbb{R}$ be functions. Given constants $A \geq 1$ and $B \geq 0$, we write $f \stackrel{A,B}{\succ} g$ to mean $f(x) \geq \frac{1}{A}g(x)-B$ for all $x \in X$, and will just write $f \succ g$ when the constants are understood.

\subsection{Curves and subsurfaces}

Throughout this paper, we let $S=S_{g,n}$ denote a connected, oriented surface of genus $g$ with $n$ punctures. Define the \textit{complexity of $S$} to be $\xi(S)=3g-3+n$. We will always assume $\xi(S)\geq 1$. Additionally,  we fix a complete hyperbolic metric on $S$. That is, we assume that $S$ is of the form $S=\mathbb{H}^2 / \Lambda$, where $\Lambda \subseteq \text{Isom}^+(\mathbb{H}^2)$ and $\Lambda$ acts properly discontinuously and freely on $\mathbb{H}^2$. 

For $i=1,2$ let $\wt{\gamma_i}$ be a bi-infinite path in $\mathbb{H}^2$ with ends limiting to distinct points $x_i$ and $y_i$ on $\partial \mathbb{H}^2$. We say that $\wt{\gamma_1}$ and $\wt{\gamma_2}$ \textit{link} if the geodesic connecting $x_1$ to $y_1$ intersects the geodesic connecting $x_2$ to $y_2$ in the interior of $\Hy$.  

By a \textit{curve} in $S$,  we will always mean the geodesic representative in the homotopy class of an essential, simple, closed curve in $S$.  By a \textit{multicurve} in $S$, we will always mean a collection of pairwise disjoint curves in $S$. We write $i(\alpha, \beta)$ to denote the geometric intersection number of curves $\alpha$ and $\beta$. 
We say that a pair of curves $\alpha$ and $\beta$ \textit{fill} $S$ if for every curve $\gamma$ in $S$ we have $i(\gamma,\alpha) > 0$ or $i(\gamma, \beta) > 0$. 

A \textit{non-annular subsurface  $Y$ of $S$} is a component of $S$ after removing a  (possibly empty) collection of pairwise disjoint  curves on $S$. Additionally, we require that $Y$ satisfies $ \xi(Y)\geq 1$; in particular, we do not consider a pair of pants to be a subsurface.  We define $\partial Y$ to be the collection of curves in $S$ that are disjoint from $Y$ and also are  contained in the closure of $Y$, treating $Y$ as a subset of $S$. When $Y \neq S$, the path metric completion of $Y$ is a surface with boundary, and the image of this boundary under the map induced by the inclusion $Y \subseteq S$ is $\partial Y$. 

An annular subsurface of $S$ is define as follows. Let $\alpha$ be a curve in $S$. Choose a component $\wt{\alpha}$ of the preimage of $\alpha$ in $\mathbb{H}^2$, and let $h \in \Lambda$ be a primitive isometry with axis $\wt{\alpha}$. Define 
\[ Y=(\overline{\mathbb{H}^2}-\{x,y\})/\langle h \rangle,\]
where $x$ and $y$ are the fixed points of $h$ on $\partial \mathbb{H}^2$. Observe that $Y$ is a compact annulus and  $\text{int}(Y)\rightarrow S$ is a covering. We say that $Y$ is the \textit{annular subsurface of $S$ with \text{core curve} $\alpha$}. We define $\partial Y$ to be $\alpha$. 

For any subsurface $Y$ of $S$, we will write $Y \subseteq S$, even though when $Y$ is an annulus, $Y$ is not a subset of $S$.

Given $f \in \Mod(S)$ and a curve or simple bi-infinite geodesic $\gamma$ in $S$, we define $f(\gamma)$ to be the curve or simple bi-infinite geodesic obtained as follows. Consider a component $\wt{\gamma}$ of the preimage of $\gamma$ in $\mathbb{H}^2$. Choose a representative $\psi$ in the isotopy class of $f$ and lift it to a map $\wt{\psi}:\mathbb{H}^2 \rightarrow \mathbb{H}^2$. We define $f(\gamma)$ to be the image in $S$ of the geodesic in $\mathbb{H}^2$ that connects the endpoints of $\wt{\psi}(\wt{\gamma})$ on $\partial \mathbb{H}^2$.  
Given $Y \subseteq S$, if $Y$ is non-annular we let $f(Y)$ denote the non-annular subsurface in its isotopy class. If $Y$ is an annulus with core curve $\alpha$, we let $f(Y)$ denote the annular subsurface of $S$ with core curve $f(\alpha)$.

\subsection{Curve complex} Let $Y$ be a subsurface of $S$.
If $Y$ satisfies $\xi(Y) \geq 1$, the \textit{curve complex of $Y$}, denoted $\mathcal{C}(Y)$, is the simplicial complex whose vertices are  curves contained in $Y$, and if $\xi(Y) >1$,  a set of vertices form a simplex if and only if they are pairwise disjoint. If $\xi(Y)=1$, then we define the simplices of $\C(Y)$ differently. In the case that $Y$ is a once punctured torus, a set of vertices form a simplex if and only if they pairwise intersect exactly once. If $Y$ is a four times punctured sphere, a set of vertices form a simplex if and only if they pairwise intersect exactly twice.

Now let $Y$ be a compact annulus. Consider all embedded arcs in $Y$ that connect one boundary component to the other. We define two arcs to be equivalent if one can be homotoped to the other, fixing the endpoints of the arcs throughout the homotopy. In this case, the \textit{curve complex of $Y$} is the simplicial complex whose vertices are equivalence classes of arcs, and a set of vertices form a simplex if and only if for each pair of vertices there exist representative arcs of each whose restrictions to $\text{int}(Y)$ are disjoint. The following simple formula will be useful to us:  given inequivalent arcs $\alpha, \beta$ in $\C(Y)$,
\begin{equation} \label{eq:annulusDistance}
d_{\C(Y)}(\alpha, \beta)=|\alpha \cdot \beta|+1, 
\end{equation}
where $\alpha \cdot \beta$ denotes the algebraic intersection number of $\alpha$ and $\beta$. 

\subsection{Markings and subsurface projection} \label{subsec:SubProjection}
A \textit{marking $\mu$ on $S$} is a maximal collection of pairwise disjoint  curves in $S$, denoted $\text{base}(\mu)$, together with another collection of associated curves called \textit{transversals}:  for each $\beta \in \text{base}(\mu)$ its associated transversal $\gamma_\beta$ is a curve that intersects $\beta$ minimally (i.e. once or twice) and is disjoint from all other curves in $\text{base}(\mu)$. 

Let $Y$ be a subsurface of $S$ and $\beta$ a multicurve in $S$. We will now define \textit{the projection of $\beta$ to $Y$}, which we will denote by $\pi_Y(\beta)$. Suppose $Y$ is not an annulus and $\beta$ is a single curve.  If $\beta$ is disjoint from $Y$, define $\pi_Y(\beta)=\emptyset$. If $\beta$ is contained in $Y$, define $\pi_Y(\beta)=\beta$. Otherwise, $\beta \cap Y$ is a collection of essential arcs in $Y$ with endpoints on $\partial Y$.  For each such arc $\gamma$, take the geodesic representatives of the boundary components of a small regular neighborhood of $\gamma \cup \partial Y$ that are contained in $Y$. Define   $\pi_Y(\beta)$ to be the collection of all such curves over all arcs $\gamma$ in $\beta \cap Y$. 
 If $\beta$ is a multicurve, define $\pi_Y(\beta)$ to be the union of the projections to $Y$ of each curve in $\beta$. 

Now let $Y$ be an annular subsurface with core curve $\alpha$ and $\text{int}(Y) \rightarrow S$ the associated covering. Let $\beta$ be a  multicurve or a bi-infinite, simple geodesic in $S$. 
 Consider the full preimage of $\beta$ in $\text{int}(Y)$. Each component is an arc in $\text{int}(Y)$ which we view as having endpoints on the boundary of $Y$.  In this case, we define $\pi_Y(\beta)$ to be the (equivalence classes of) arcs in this collection that have an endpoint on each boundary component of $Y$.  When convenient, we will write $\pi_\alpha(\beta)$ instead of $\pi_Y(\beta)$. 

We now describe how to project a marking $\mu$ to $Y\subseteq S$. If $Y$ is non-annular or $Y$ is an annulus whose core curve is not contained in $\text{base}({\mu})$, we define $\pi_Y(\mu)=\pi_Y(\text{base}(\mu))$. Otherwise, $Y$ is an annulus with core curve $\alpha \in \text{base}(\mu)$, and we define $\pi_Y(\mu)$ to be $\pi_Y(\gamma_\alpha)$, where $\gamma_\alpha$ is the transversal associated to $\alpha$. 

Given any subsurface $Y \subseteq S$, we define \[d_Y(\mu, \mu')= \text{diam}_{\mathcal{C}(Y)}(\pi_Y(\mu) \cup \pi_Y(\mu')),\] where $\mu$ and $\mu'$ are markings, collections of curves, or (when $Y$ is an annulus) bi-infinite simple geodesics in $S$. A useful fact about subsurface projection is the following.  For all $f\in \Mod(S)$ 
\[d_{Y}(f(\mu), f(\mu'))=d_{f^{-1}(Y)}(\mu, \mu').\]

In this paper, we utilize the following theorem, which involves subsurface projections.

\begin{theorem} [Lemma 2.3 in \cite{MMII}] \label{thrm:boundProjmultiCurve}
For all subsurfaces $W$ of $S$, given any marking or multicurve $\mu$ such that $\pi_W(\mu) \neq \emptyset$, we have that $\text{diam}_{\C(W)}(\pi_W(\mu)) \leq 2$.  If $W$ is an annulus, then $\text{diam}_{\C(W)}(\pi_W(\mu)) \leq 1$.
\end{theorem} 

Masur and Minsky \cite{MMII} define the \textit{marking graph of $S$}, denoted $\wt{\mathcal{M}}(S)$, to be the graph whose vertices are markings and vertices are adjacent if one can be obtained from the other by an elementary move; see \cite{MMII} for a complete definition. Giving $\wt{\mathcal{M}}(S)$ the path metric $d_{\wt{\mathcal{M}}(S)}$ and $\Mod(S)$ a word metric $d_{\Mod(S)}$,  there is an action of $\Mod(S)$ on $\wt{\mathcal{M}}(S)$ by isometries for which every orbit map is a quasi-isometry. The following theorem gives a relationship between distances in $\wt{\mathcal{M}}(S)$ and subsurface projections. 

\begin{theorem}[Lemma 3.5 in \cite{MMII}] \label{theorem:MMLipProjection}
For any subsurface $W$ of $S$ and any markings $\mu$ and $\mu'$ on $S$, we have that $d_W(\mu, \mu') \leq 4 d_{\wt{\mathcal{M}}(S)}(\mu, \mu')$. 
\end{theorem}

We say that distinct subsurfaces $X$ and $Y$ are \textit{disjoint} if $\pi_X(\partial Y)=\emptyset$ and $\pi_Y(\partial X)=\emptyset$. We say that $X$ is a \textit{proper subsurface of $Y$}, denoted $X \subsetneq Y$, if $\pi_Y(\partial X)\neq \emptyset$ and $\pi_X(\partial Y) =\emptyset$. We say that $X$ and $Y$ are \textit{overlapping}, denoted $X \pf Y$,  if $\pi_Y(\partial X) \neq \emptyset$ and $\pi_X(\partial Y) \neq \emptyset$. In the case where $X$ and $Y$ are not annuli,  these relationships, respectively, are disjointness, proper containment, and intersection without containment as subsets of $S$.
We say  $X$ and $Y$ \textit{fill} $S$ if for every curve $\gamma$ in $S$ we have $\pi_X(\gamma) \neq \emptyset$ or $\pi_Y(\gamma) \neq \emptyset$. 

The following theorems will be used to prove our results. The first theorem was proved in \cite{Behrstock} and later a simpler proof with constructive constants appeared in \cite{Mangahas13}. 
\begin{theorem}
[Behrstock inequality: Theorem 4.3 in  \cite{Behrstock}, Lemma 2.13 in \cite{Mangahas13}]
\label{bineq}
Let $X$ and $Y$ be overlapping subsurfaces of $S$ and $\mu$ a marking on $S$. Then
\[ d_X(\mu, \partial Y) \geq 10 \hspace{10pt} \text{ implies that } \hspace{10pt} d_Y(\mu, \partial X) \leq 4.\]
\end{theorem}

\begin{theorem} [Bounded Geodesic Image Theorem: Theorem 3.1 in \cite{MMII}] \label{bgit} There exists a constant $K_0$ depending only on $S$ such that the following is true. 
Let $Y$ and $Z$ be subsurfaces of $S$ with $Y$ a proper subsurface of $Z$. Let $v_1, \ldots, v_n$ be any geodesic segment in $\mathcal{C}(Z)$ satisfying $\pi_Y(v_i) \neq \emptyset$ for all $1\leq i \leq n$. Then 
\[\text{diam}_{\C(Y)}(\pi_Y(v_1) \cup \ldots \cup \pi_Y(v_n)) \leq K_0.\]
\end{theorem} 

\subsection{Partial order on subsurfaces} Let $\mu, \mu'$ be markings on $S$ and $K \geq 20$. Let $\Omega(K, \mu, \mu')$ denote the collection of subsurfaces $Y$ of $S$ such that $d_Y(\mu, \mu') \geq K$. 
Behrstock, Kleiner, Minsky, and Mosher \cite{BKMM} define the following partial order on $\Omega(K, \mu, \mu')$.
Given $X, Y \in \Omega(K, \mu, \mu')$ such that $X \pf Y$,  define $X \prec Y$ if and only if one of the following equivalent conditions is satisfied:
\[ d_X(\mu, \partial Y) \geq 10, \hspace{9pt} d_X(\partial Y, \mu') \leq 4, \hspace{9pt} d_Y(\mu, \partial X) \leq 4, \hspace{9pt} \text{ or } \hspace{2pt} d_Y(\partial X, \mu') \geq 10.\] That these conditions are equivalent is a consequence of Theorem \ref{bineq}; see Corollary 3.7 in \cite{CLM}.

\subsection{Embedding RAAGs in Mod(S)}
If $f \in \Mod(S)$ is such that there exists a representative in the isotopy class of $f$ that pointwise fixes the complement of a non-annular subsurface $Y$, we say that $f$ is \textit{supported on $Y$}.
Given such an $f$,  we define the \textit{translation length of $f$ on $\C(Y)$} to be
\[\tau_Y(f)=\lim_{n\rightarrow \infty} \frac{d_Y(\mu, f^n(\mu))}{n},\]
where $\mu$ is any marking on $S$. If $f\in \Mod(S)$ is a power of a Dehn twist about a curve $\alpha$, we say that $f$ is \textit{supported on} the annular subsurface $Y$ with core curve $\alpha$, and define $\tau_Y(f)$ to be the absolute value of the power.
In either case, we say that $Y$ \textit{fully supports} $f$ if $\tau_Y(f) >0.$ By the work of Masur and Minsky \cite{MMI}, when $Y$ is non-annular, $Y$ fully supports $f$ if and only if  $f$ is pseudo-Anosov on $Y$.

 Clay, Leininger, and Mangahas \cite{CLM} proved the following result, which allows us to find quasi-isometrically embedded right-angled Artin subgroups inside $\Mod(S)$.  

\begin{theorem} [Theorem 2.2 in \cite{CLM}]  \label{CLMemb} Let $\Gamma$ be a finite graph with $V(\Gamma)=\{s_1, \ldots, s_k\}$, and let  $\{X_1, \ldots, X_k\}$ be a collection of non-annular subsurfaces of $S$. Suppose $s_is_j$ is an edge in $\Gamma$ if and only if $X_i$ and $X_j$ are disjoint, and $s_is_j$ is not an edge in $\Gamma$ if and only if $X_i \pf X_j$ or $i =j$. 
Then there exists a  constant $C>0$ such that the following holds. Let $\{f_1, \ldots, f_k\}$ be a set of mapping classes of $S$ such that $f_i$ is  pseudo-Anosov on $X_i$ and satisfies $\tau_{X_i}(f_i) \geq C$ for all $i$. Then the homomorphism \[\phi: A(\Gamma) \rightarrow \Mod(S) \hspace{10pt} \text{ defined by } \hspace{10pt} \phi(s_i)=  f_i \hspace{10pt} \text{ for all } \hspace{3pt}i\] is a quasi-isometric embedding, and in particular is injective.
\end{theorem}

Koberda \cite{Koberda} also has a result which produces right-angled Artin subgroups of $\Mod(S)$. Below we give a special case of Koberda's result that we will use.

\begin{theorem}[Theorem 1.1 in \cite{Koberda}] \label{Kemb}
Let $\{\alpha_1, \ldots, \alpha_k\}$ be a collection of distinct curves in $S$. Let $\Gamma$ be the graph with $V(\Gamma)=\{s_1, \ldots, s_k\}$ and with $s_is_j$ an edge in $\Gamma$ if and only if  $i(\alpha_i, \alpha_j)=0$. Then for sufficiently large $N$, the homomorphism 
\[\phi: A(\Gamma) \rightarrow \Mod(S)\hspace{10pt} \text { defined by } \hspace{10pt}  \phi(s_i)=  T^N_{\alpha_i}  \hspace{10pt} \text{ for all } \hspace{3pt} i,\]
is injective, where $T_{\alpha_i}$ denotes a Dehn twist about $\alpha_i$. 
\end{theorem}
\subsection{Gromov boundary of hyperbolic spaces} \label{subsec:GromovHyp}

A geodesic metric space $X$ is \textit{Gromov hyperbolic} (or just hyperbolic) if there exists a $\delta \geq 0$ such that given any geodesic triangle in $X$, each side is contained in the $\delta$-neighborhood of the union of the other two sides. Given a Gromov hyperbolic space $(X,d_X)$ and points $x,y,z \in X$, the \textit{Gromov product} of $x$ and $y$ with respect to $z$ is defined as

\[(x,y)_z=\frac{1}{2}\left(d_X(x,z)+d_X(y,z)-d_X(x,y)\right).\]

We say that a sequence $(x_n)$ in $X$ \textit{converges at infinity} if $\displaystyle{\liminf_{i,j \rightarrow \infty} (x_i,x_j)_z}=\infty$ for some (any) $z \in X$. We define two such sequences $(x_n)$ and $(y_n)$ to be equivalent if $\displaystyle{\liminf_{i,j \rightarrow \infty} (x_i, y_j)_z =\infty}$ for some (any) $z \in X$. 
The \textit{Gromov boundary of $X$} is the collection of all such sequences up to this equivalence, and is denoted $\partial_GX$ or just $\partial X$ when it is clear from context that we are using the Gromov boundary. 

One Gromov hyperbolic space that this paper is concerned with is the curve complex of $S$, which was proved to be Gromov hyperbolic by Masur and Minsky \cite{MMII}. We can now state a corollary of Theorem \ref{bgit} that will be useful later.  
\begin{corollary} \label{bgitCor}
Let $X$ and $Y$ be subsurfaces of $S$ with $X$ a proper subsurface of $Y$. Suppose $(\mu_n)_{n \in \mathbb{N}}$ is a sequence of markings on $S$ such that $\pi_Y(\mu_n) \rightarrow \lambda$ for some $\lambda \in \partial \mathcal{C}(Y)$.  Then $\text{diam}_{\C(X)}(\pi_X(\mu_1) \cup \pi_X(\mu_2) \cup  \ldots ) < \infty.$
\end{corollary}

\begin{proof}
For each $n$, choose  $\alpha_n \in \pi_Y(\mu_n)$. Because $\pi_Y(\mu_n) \rightarrow \lambda \in \partial \mathcal{C}(Y)$, we can choose $L$ large so that for all $n \geq L$ we have
\begin{equation} \label{eq:GromovProdBound}
(\alpha_n, \alpha_L)_{\alpha_1} \geq 2+ d_Y(\partial X, \alpha_1),
\end{equation} 
where the Gromov product is computed in $\C(Y)$.  
Consider $n \geq L$. Let $\gamma_n$ be a geodesic in $\C(Y)$ with endpoints $\alpha_n$ and $\alpha_L$. If there exists a vertex $v$ on $\gamma_n$ with $\pi_X(v) = \emptyset$, then $v$ and $\partial X$ form a multicurve in $Y,$ which implies that
\begin{align*}(\alpha_n, \alpha_L)_{\alpha_1} 
&=\frac{1}{2}\bigg(d_Y(\alpha_n,\alpha_1)+d_Y(\alpha_L, \alpha_1)-d_Y(\alpha_n, \alpha_L)\bigg) \\ 
& \leq \frac{1}{2} \bigg( d_Y(\alpha_n,v)+d_Y(v, \alpha_1) +d_Y(\alpha_L, v)+d_Y(v, \alpha_1)-(d_Y(\alpha_n, v)+d_Y(v, \alpha_L)) \bigg)\\
&=d_Y(v, \alpha_1) \leq d_Y(v,\partial X)+d_Y(\partial X, \alpha_1) \leq 1+d_Y(\partial X, \alpha_1).
\end{align*}
But this contradicts Inequality (\ref{eq:GromovProdBound}), so we conclude that $\pi_X(v) \neq \emptyset$ for all $v$ on $\gamma_n$. We can now apply Theorems \ref{thrm:boundProjmultiCurve} and  \ref{bgit} to see that for all $n \geq L$
\[d_X(\mu_n, \mu_L) \leq \text{diam}_{\C(X)}(\pi_X(\mu_n))+d_X(\alpha_n, \alpha_L)+\text{diam}_{\C(X)}(\pi_X(\mu_L)) \leq 2+K_0+2,\]
where $K_0$ is an in Theorem \ref{bgit}. Therefore, 
\[\text{diam}_{\C(X)}(\pi_X(\mu_1)\cup \pi_X(\mu_2)\cup \ldots ) \leq \text{diam}_{\C(X)}( \pi_X(\mu_1)\cup \ldots \cup \pi_X(\mu_L))+2(K_0+4)<\infty.\]
\end{proof}

\subsection{Hierarchically hyperbolic spaces} \label{subsec:HHS}

In \cite{BHS} Behrstock, Hagen, and Sisto define the notion of a hierarchically hyperbolic space. Roughly, a hierarchically hyperbolic space is a quasi-geodesic metric space $\mathcal{X}$, equipped with additional structure which we will call a hierarchically hyperbolic space (HHS) structure. An HHS structure consists of an index set $\mathcal{G}$ and for each $W \in \mathcal{G}$ a Gromov hyperbolic space $\widehat{C}W$ and a projection map $\pi_W: \mathcal{X} \rightarrow2^{\widehat{C}W}$. The elements of $\mathcal{G}$ and the projection maps must satisfy a long list of properties. See \cite{BHS} and \cite{BHSII}.

 The first example of a hierarchically hyperbolic space is $\Mod(S)$, where here $\mathcal{G}$ is the collection of all subsurfaces of $S$, $\widehat{C}W$ is the curve graph of $W$ for $W \in \mathcal{G}$, and projection $\pi_W$ is given by composing an orbit map for the action of $\Mod(S)$ on $\wt{\mathcal{M}}(S)$ with the subsurface projection map defined in Section \ref{subsec:SubProjection}. The work of Masur and Minsky  \cite{MMI},\cite{MMII} and Behrstock \cite{Behrstock}  imply that $\Mod(S)$ is a hierarchically hyperbolic space. See \cite{BHSII} Section 11 for details. In fact, the notion of hierarchical hyperbolicity was motivated by a desire to generalize some of the machinery surrounding mapping class groups. 
 
In \cite{BHS} it is shown that a large class of $\text{CAT}(0)$ cube complexes can be equipped with a hierarchically hyperbolic structure, including the universal covers of Salvetti complexes associated to right-angled Artin groups. The $\text{CAT}(0)$ cube complex we are primarily concerned with is the Cayley graph $X$ of $A(\Gamma)$ when $\Gamma$ has no edges (that is, $A(\Gamma)$ is a free group). 
We equip $A(\Gamma)$ with a hierarchically hyperbolic structure by equipping $X$ with such a structure and then associating $A(\Gamma)$ with $X$.

\subsection{Boundary of hierarchically hyperbolic spaces} \label{subsec:BoundaryHHS}
In \cite{DHS} the authors construct a boundary for hierarchically hyperbolic spaces. Here we will describe convergence in this boundary for $\Mod(S)$ and for free groups. With the exception of Theorem \ref{theorem:NofillNoextend}, these will be the only examples we will need.

As a set, the HHS boundary of $\Mod(S)$ is defined as follows:
\begin{multline}
\partial \Mod(S) =\Bigg\{\sum_{Y \subseteq S} c_Y\lambda_Y : c_Y \geq 0 \text{ and } \lambda_Y \in \partial \mathcal{C}(Y) \text{ for all } Y,   \sum_{Y \subseteq S} c_Y =1, \Bigg. \\ \Bigg. \text{ and if } c_{Y'}, c_{Y}>0, \text{ then } Y\text{ and } Y' \text{ are disjoint or equal} \Bigg\}. \nonumber
\end{multline}
 In \cite{DHS}, the authors define a topology on $\Mod(S) \cup \partial \Mod(S)$. In this topology, Definition 2.10 of \cite{DHS} tells us that a sequence of mapping classes $(g_n)_{n \in \mathbb{N}}$ in $\Mod(S)$ converges to a point $\displaystyle{\sum_{i=1}^k c_i \lambda_i}$ in $\partial \Mod(S)$, where $c_i >0$ for all $i$,  $\displaystyle{\sum_{i=1}^k c_i =1}$, and $\lambda_i \in \partial \mathcal{C}(Y_i)$ for pairwise disjoint subsurfaces $Y_1, \ldots, Y_k$, if and only if the following statements hold: For a fixed marking $\mu$ on $S$, 

\begin{enumerate}
\item $\displaystyle{\lim_{n \rightarrow \infty} \pi_{Y_i}(g_n \mu) = \lambda_i}$ for each $i=1, \ldots, k$,
\item $\displaystyle{\lim_{n \rightarrow \infty} \frac{d_{Y_i}(\mu, g_n \mu)}{d_{Y_j}(\mu, g_n \mu)}=\frac{c_i}{c_j}}$ for each $i,j=1, \ldots, k$, and
\item $\displaystyle{\lim_{n\rightarrow \infty} \frac{d_W(\mu, g_n \mu)}{d_{Y_i}(\mu, g_n \mu)}=0}$ for every (any) $i=1, \ldots, k$ and every subsurface $W \subseteq S$ that is disjoint from $Y_j$ for all $j=1, \ldots, k$. 
\end{enumerate}

Let $\Gamma$ be a graph with no edges, and let $A(\Gamma)$ be the corresponding free group, equipped with an HHS structure. The HHS boundary of $A(\Gamma)$ will be denoted by $\partial A(\Gamma)$. We do not define $\partial A(\Gamma)$ here because Theorem 4.3 in \cite{DHS} implies that the identity map $A(\Gamma) \rightarrow A(\Gamma)$ extends to a homeomorphism $A(\Gamma) \cup \partial_GA(\Gamma) \rightarrow A(\Gamma) \cup \partial A(\Gamma)$. Thus, two sequences in $A(\Gamma)$ converge to the same point in $\partial_G A(\Gamma)$ if and only if they converge to the same point in $\partial A(\Gamma)$. (See Section 2 of \cite{DHS} for the definition of $\partial A(\Gamma)$.)

Another useful fact on convergence is that $\Mod(S) \cup \partial \Mod(S)$ and $A(\Gamma) \cup \partial A(\Gamma)$ are sequentially compact (see Theorem 3.4 of \cite{DHS}). 

To understand Question \ref{mainq} and the statements of our theorems, one last definition is needed. 

\begin{definition} \normalfont
Let $\phi: A(\Gamma) \rightarrow \Mod(S)$ be an injective homomorphism and let $A(\Gamma)$ and $\Mod(S)$ be equipped with any fixed HHS structures. We say that $\phi$ \textit{extends continuously to a map} $\partial A(\Gamma) \rightarrow \partial \Mod(S)$ if there exists a function $\overline{\phi}: A(\Gamma) \cup \partial A(\Gamma) \rightarrow \Mod(S) \cup \partial \Mod(S)$ such that 
(1) $\overline{\phi}|_{A(\Gamma)}=\phi$, 
(2) $\overline{\phi}(\partial A(\Gamma)) \subseteq \partial \Mod(S)$, and (3) $\overline{\phi}$ is continuous at each point in $\partial A(\Gamma)$. 
\end{definition}

\begin{remark}
\normalfont  To establish that  $\phi: A(\Gamma) \rightarrow \Mod(S)$ extends continuously, it is enough to show that for all $x \in \partial A(\Gamma)$, given any two sequences $(x_n)$ and $(y_n)$ in $A(\Gamma)$ that converge to $x$, we have that $(\phi(x_n))$ and $(\phi(y_n))$ converge to the same point in $\partial \Mod(S)$. This follows from a diagonal sequence argument (see the end of the proof of Theorem 5.6 in \cite{DHS} for details).
\end{remark}

\section{Lemmas on subsurface projections} \label{sec:lemma}
The following lemmas are the heart of our proof of Theorem \ref{theoremCLMKintro}. 

\begin{lemma} \label{lemann}
 Suppose $X$ and  $Y$ are disjoint subsurfaces of $S$, and if $Y$ is an annulus, then the core of $Y$ is not contained in $\partial X$.  If $\mu$ and $\mu'$ are markings and  $f\in Mod(S)$ a mapping class supported on $X$, then $|d_Y(\mu, f(\mu')) - d_Y(\mu, \mu')|  \leq 4$. 
\end{lemma}

\begin{proof}
If $Y$ is not an annulus, then $\pi_Y(f(\mu'))=\pi_Y(\mu')$ so the claim clearly holds. Assume then that $Y$  is an annular subsurface of $S$ with core $\alpha$, and let $\text{int}(Y) \rightarrow S$ be the associated covering.  
Because $\pi_X(\alpha)=\emptyset$ and $\alpha$ is not in $\partial X$, we can find a curve $\gamma$ in $S$, distinct from $\alpha$, that intersects $\alpha$ and satisfies $\pi_X(\gamma)=\emptyset$. 
If $X$ is not an annulus, define $Z$ to be the component of $S-X$ that contains $\alpha$. If $X$ is an annulus with core $\beta$, let $Z$ be the component of $S$ containing $\alpha$ after removing a small regular neighborhood of $\beta$. The neighborhood should be taken small enough so that $\gamma$ is contained in $Z$. Let $\wt{\alpha}$ be the component of the preimage of $\alpha$ in $\text{int}(Y)$ that is a closed curve. 
Let $\wt{Z}$ be the component of the preimage of $Z$ in $\text{int}(Y)$ that contains $\wt{\alpha}$.

Abusing notation, we let $f$ denote a representative in the isotopy class of $f$ that fixes $Z$ pointwise. 
 Let $\wt{f}: \text{int}(Y) \rightarrow \text{int}(Y)$ denote the lift of  $f$  that fixes a point on $\wt{\alpha}$, and thus fixes $\wt{Z}$ pointwise.  Let $\wt{\gamma}$ be a component of the preimage of $\gamma$ in $\text{int}(Y)$ that  intersects $\wt{\alpha}$. Then $\wt{\gamma}$ is contained in $\wt{Z}$, implying that $\wt{f}$ fixes $\wt{\gamma}$ pointwise.
 
It then follows that for $\beta' \in \pi_Y(\mu')$, we have (after replacing $\wt{f}(\beta')$ with a some representative in its homotopy class) that $\wt{f}(\beta')$ and $\beta'$ intersect at most once.  Thus, by Equation (\ref{eq:annulusDistance}) we have $d_{\mathcal{C}(Y)}(\wt{f}(\beta'),\beta) =1 +|\wt{f}(\beta') \cdot \beta'| \leq 2$. 
Now apply the triangle inequality and Theorem \ref{thrm:boundProjmultiCurve} to see that 
\begin{align*} |d_Y(\mu, f(\mu')) -d_Y(\mu, \mu')| 
& \leq d_Y(\mu', f(\mu')) \\
&\leq \text{diam}_{\C(Y)}(\pi_Y(\mu')) +d_{C(Y)}(\beta', \wt{f}(\beta'))+\text{diam}_{\C(Y)}(\pi_Y(f(\mu'))) \\
&\leq 1 +2+1 =4.
\end{align*}
\end{proof}

 \begin{lemma} \label{lemmaWordLength} Given a homomorphism $\phi: A(\Gamma) \rightarrow \Mod(S)$ and a marking $\mu$ on $S$, there exists a constant $M \geq 1$ such that the following holds. 
 Let $y_1\ldots y_n \in A(\Gamma)$, where each $y_i \in V(\Gamma)$. Then $d_W(\mu, \phi(y_1\ldots y_n) \mu) \leq Mn$ for all subsurfaces $W \subseteq S$. 
 \end{lemma}

\begin{proof}
 Define $\displaystyle{M=4\max\{ d_{\wt{\mathcal{M}}(S)}(\mu, \phi(x)\mu): x \in V(\Gamma)\}}$. By the triangle inequality and Theorem \ref{theorem:MMLipProjection}, 
 \begin{align*} 
 d_W(\mu, \phi(y_1\ldots y_n) \mu) 
 & \leq \sum_{i=1}^nd_W(\phi(y_1 \ldots y_{i-1})\mu, \phi(y_1 \ldots y_i)\mu) \\
 &  \leq  \sum_{i=1}^n 4d_{\wt{\mathcal{M}}(S)}(\phi(y_1 \ldots y_{i-1})\mu, \phi(y_1 \ldots y_i)\mu) \\
& = \sum_{i=1}^n 4d_{\wt{\mathcal{M}}(S)}(\mu, \phi(y_i)\mu)   \leq Mn.
 \end{align*}  
\end{proof}

\begin{lemma} \label{lemmaGeneralLin}
Let $\phi: A(\Gamma) \rightarrow \Mod(S)$ be a homomorphism. Let $(g_n)_{n \in \mathbb{N}}$ be a sequence of elements in $A(\Gamma)$ and $\mu$ a  marking on $S$. Suppose for some subsurface $W \subseteq S$ there exist constants $A\geq 1$ and $B \geq 0$, that do not depend on $n$, such that $d_W(\mu, \phi(g_n)\mu) \stackrel{A,B}{\succ} ||g_n||,$ where $|| g_n||$ denotes the word length of $g_n$ with respect to the standard generating set $V(\Gamma)$ for $A(\Gamma)$. Further suppose that $\displaystyle{\lim_{n \rightarrow \infty}||g_n|| = \infty}$ and that $(\pi_W(\phi(g_n) \mu))_{n \in \mathbb{N}}$ converges to some point $\lambda_W$ in $\partial \mathcal{C}(W)$. 
Then all accumulation points of $(\phi(g_n))_{n \in \mathbb{N}}$ in $\Mod(S) \cup \partial \Mod(S)$  are in $\partial \Mod(S)$ and are of the form $\displaystyle{\sum_{Y \subseteq S} c_Y \lambda_Y}$, where $c_W > 0$.  
\end{lemma}

\begin{proof}  After passing to a subsequence, we may assume that $(\phi(g_n))_{n \in \mathbb{N}}$ converges.  
 By \linebreak  assumption, 
$\displaystyle{\lim_{n \rightarrow \infty} d_W(\mu, \phi(g_{n}) \mu) =\infty}$. Combine this with Theorem \ref{theorem:MMLipProjection} to see that \linebreak  $\displaystyle{\lim_{n \rightarrow \infty} d_{\widetilde{\mathcal{M}}(S)}( \mu, \phi(g_{n}) \mu) =\infty}$. Because $\wt{\mathcal{M}}(S)$ is quasi-isometric to $\Mod(S)$ via orbit maps,  it follows that  $\displaystyle{\lim_{n \rightarrow \infty} d_{\Mod(S)}( 1, \phi(g_{n})) =\infty}$.  Thus, it must be that $\displaystyle{\lim_{n \rightarrow \infty} \phi(g_{n}) \in \partial \Mod(S)}$. 

Suppose $\displaystyle{\lim_{n\rightarrow \infty} \phi(g_{n})=\sum_{Y \subseteq S} c_Y \lambda_Y}$ for constants $c_Y \geq 0$ and $\lambda_Y \in \partial \mathcal{C}(Y)$. We will now argue that $c_W >0$.  Let $Z \subseteq S$ be such that $c_Z >0$. If $W=Z$, we are done. So we assume $W \neq Z$. By definition of the topology on $\Mod(S) \cup \partial \Mod(S)$, we have that $\displaystyle{\lim_{n \rightarrow \infty} \pi_Z( \phi(g_{n})\mu)=\lambda_Z}$. If $W \subsetneq Z$, then Corollary \ref{bgitCor} implies that $ \text{diam}_{\C(W)}(\pi_W(\phi(g_{1})\mu) \cup \pi_W(\phi(g_{2})\mu) \cup \ldots) < \infty$. But this cannot be since $\pi_W(\phi(g_{n})\mu) \rightarrow \lambda_W \in \partial \C(W)$. Similarly, we cannot have $Z \subsetneq W$ for then Corollary \ref{bgitCor} implies that  $\text{diam}_{\C(Z)}( \pi_Z(\phi(g_{1})\mu) \cup \pi_Z(\phi(g_{2})\mu) \cup\ldots ) < \infty$, contradicting that  $\pi_Z( \phi(g_{n})\mu) \rightarrow \lambda_Z \in \partial \C(Z)$. Now suppose that $Z \pf W$. Then by Theorem \ref{bineq}, after passing to a subsequence,  we have that
\[d_W(\partial Z, \phi(g_{n}) \mu) \leq 10 \hspace{5pt} \text{ for all } n, \hspace{10pt} \text{ or } \hspace{10pt} d_Z(\partial W, \phi(g_{n}) \mu) \leq 10 \hspace{5pt} \text{ for all } n.\]
If $d_W(\partial Z, \phi(g_{n}) \mu) \leq 10$ for all $n$, then for all $n$
\[d_W(\mu, \phi(g_{n}) \mu) \leq d_W(\mu, \partial Z) +d_W(\partial Z, \phi(g_{n}) \mu) \leq d_W(\mu, \partial Z) +  10,\]
contradicting that $\pi_W(\phi(g_{n})\mu)  \rightarrow \lambda_W \in \partial \C(W)$. Similarly, if $d_Z(\partial W, \phi(g_{n}) \mu) \leq 10$ for all $n$, then $d_Z(\mu, \phi(g_{n}) \mu)$ is bounded independent of $n$ contradicting that $\pi_Z(\phi(g_{n}) \mu) \rightarrow \lambda_Z \in \partial \C(Z)$. So it is not the case that $Z \pf W$. 
Therefore it must be that $W$ and $Z$ are disjoint for all $Z \subseteq S$ with $c_Z >0$. 

 Fix $Z \subseteq S$ with $c_Z >0$.
 Lemma \ref{lemmaWordLength} together with the fact that $d_W(\mu, \phi(g_{n})\mu) \stackrel{A,B}{\succ} ||g_{n}||$ implies that 
\begin{equation} \label{EqRatio}\frac{d_W(\mu, \phi(g_{n})\mu)}{d_Z(\mu, \phi(g_{n})\mu)} \geq \frac{\frac{1}{A}||g_{n}||-B}{M||g_{n}||},
\end{equation}
where $M \geq 1 $ is as in Lemma \ref{lemmaWordLength}.
Since $||g_{n}|| \rightarrow \infty$, Equation (\ref{EqRatio}) implies 
\[\lim_{n \rightarrow \infty} \frac{d_W(\mu, \phi(g_{n})\mu)}{d_Z(\mu, \phi(g_{n})\mu)} \geq \lim_{n \rightarrow \infty} \frac{\frac{1}{A}||g_{n}||-B}{M||g_{n}||} >0.\]
Therefore by definition of the topology of $\Mod(S) \cup \partial \Mod(S)$, we have $c_W >0$ as desired.

\end{proof}
 \section{Clay, Leininger, Mangahas RAAGs} \label{sec:CLMRAAG}
 
 In this section, we prove the first part of Theorem \ref{theoremCLMKintro}.  We begin with a description of a Clay, Leininger, Mangahas embedding $\phi: A(\Gamma)\rightarrow \Mod(S)$. 
  \begin{figure} 
 \center
 \def\svgwidth{275pt}
 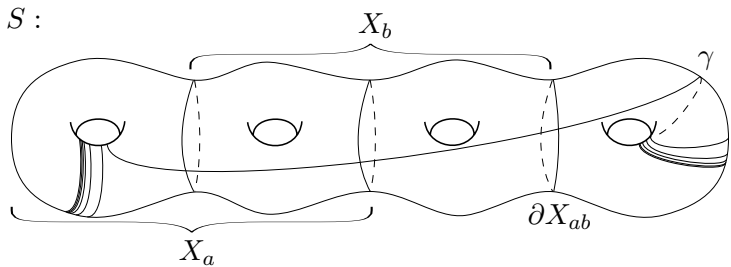
\vspace{-20pt}
\captionsetup{width=.8\linewidth}
 \caption{Overlapping subsurfaces $X_a$ and $X_b$ of surface $S$, curve $\partial X_{ab}$, and bi-infinite simple  geodesic $\gamma$.}
 \label{surface} 
 \end{figure}
 
 \textbf{Embedding construction:} Let $\Gamma$ be the graph with vertex set $V(\Gamma)=\{a, b\}$ and no edges.  Let $S=\mathbb{H}^2/\Lambda, X_a,$ and $X_b$ be the surfaces indicated in Figure \ref{surface}. 
 For short, let $X_{ab}$ denote $X_a \cup X_b$.  Let $\widetilde{S-X_{ab}}$ be a component of the preimage of $S-X_{ab}$ in $\mathbb{H}^2$, and let $\wt{\partial X_{ab}}$ be a geodesic in $\Hy$ that is in the boundary of $\wt{S-X_{ab}}$. 
 
Let $\wt{\gamma}$ be a geodesic in $\mathbb{H}^2$ that links with $\wt{\partial X_{ab}}$ and maps to a simple bi-infinite geodesic $\gamma$ in $S$. Further suppose that $\wt{\gamma} \cap (\wt{S-X_{ab}})$ is an infinite ray and let $p$ be its endpoint on $\partial \mathbb{H}^2$. For example, take $\gamma$ to be the simple bi-infinite geodesic in $S$ with one end spiraling around a curve essential in $S-X_{ab}$ and the other end spiraling around a curve in $X_a$ as in Figure \ref{surface}, and take $\wt{\gamma}$ to be an appropriate lift of $\gamma$. 
 Choose $f_b\in \Mod(S)$ so that $f_b$ is pseudo-Anosov  on $X_b$. To simplify arguments, we abuse notation and let $f_b$ denote a representative in the isotopy class of $f_b$ that fixes all points outside $X_b$. This ensures that $\tilde{f_b}$ fixes $\widetilde{S-X_{ab}}$ pointwise, where $\tilde{f_b}:\mathbb{H}^2\rightarrow \mathbb{H}^2$ is the lift of $f_b$ fixing some point on $\wt{\partial X_{ab}}$. Thus, the extension of $\wt{f_b}$ to $\partial \mathbb{H}^2$  fixes pointwise $p$ and the endpoints $x$ and $y$ of $\wt{\partial X_{ab}}$.
 Additionally, we choose $f_b$ to have the following properties:
 \begin{enumerate}
 \item $\tilde{f_b}(\wt{\gamma})$ links with $h(\wt{\gamma})$, where $h\in \Lambda$ is a primitive isometry  with axis $\widetilde{\partial X_{ab}}$, and
 \item $\tau_{X_b}(f_b) \geq C,$ where $C$ is as in Theorem \ref{CLMemb}.  
 \end{enumerate}  
 
We note that a pseudo-Anosov on $X_b$ satisfying (1) can be obtained from any mapping class that is pseudo-Anosov on $X_b$ by post-composing with some number of Dehn twists (or inverse Dehn twists) about $\partial X_{ab}$. Finally, a pseudo-Anosov on $X_b$ satisfying (1) and (2) can be obtained from one satisfying (1) by passing to a sufficiently high power. 
 
 Let $f_a \in \Mod(S)$ be any mapping class that is pseudo-Anosov on $X_a$ and satisfies $\tau_{X_a}(f_a) \geq C$. Theorem \ref{CLMemb} says that the homomorphism 
  \[\phi: A(\Gamma) \rightarrow \Mod(S)\hspace{10pt}  \text{ defined by } \hspace{10pt} \phi(a)=f_a, \hspace{5pt}  \phi(b)= f_b \] 
  is a quasi-isometric embedding.  
  
 Equip $A(\Gamma)$ with any HHS structure.
 In the remainder of this section we will prove the following theorem, which proves the first part of Theorem \ref{theoremCLMKintro}. 

\begin{theorem} \label{thrmmain}
The sequences $(a^n)_{n \in \mathbb{N}}$ and $(a^nb^n)_{n \in \mathbb{N}}$ converge to the same point in $\partial A(\Gamma)$, but $(\phi(a^n))_{n \in \mathbb{N}}$ and $(\phi(a^nb^n))_{n \in \mathbb{N}}$ do not converge to the same point in $\partial \Mod(S)$. 
\end{theorem}

We will divide the proof of Theorem \ref{thrmmain} into two propositions. 

\begin{proposition} \label{prop:CLMseqcon}
The sequences $(a^n)_{n \in \mathbb{N}}$ and $(a^nb^n)_{n \in \mathbb{N}}$ converge to the same point in  $\partial A(\Gamma)$.
\end{proposition}

\begin{proof}
Let $X$ be the Cayley graph of $A(\Gamma)$. 
By the discussion in Section \ref{subsec:BoundaryHHS}, to show that $(a^n)_{n \in \mathbb{N}}$ and $(a^nb^n)_{n \in \mathbb{N}}$ converge to the same point in  $\partial A(\Gamma)$, it is enough to show that they  converge to the same point in $\partial_GX$. Now the  Gromov product 
\[(a^i, a^j b^j)_1= \min(i,j) \rightarrow \infty \text{ as } i, j \rightarrow \infty. \]
Thus, $\displaystyle{\lim_{n \rightarrow \infty} a^n=\lim_{n\rightarrow \infty} a^n b^n}$ in $\partial_G X$, as desired. 

\end{proof}

Throughout the rest of this section $\mu$ will denote a fixed marking on $S$. To continue, we require the following lemma. 

\begin{lemma}  \label{lemlingrowth} There exist constants $A\geq 1$ and $B \geq 0$ such that for all $n \geq 1$ we have  $d_{\partial X_{ab}}(\mu, \phi(a^nb^n)\mu) \stackrel{A,B}{\succ} n$. Consequently, after passing to a subsequence, $(\pi_{\partial X_{ab}} (\phi(a^nb^n) \mu))_{n \in \mathbb{N}}$ converges to a point in $\partial \mathcal{C}(\partial  X_{ab})$. 
\end{lemma} 

\begin{proof} 

We begin by establishing the following claim. 

\

\noindent \textit{Claim 1: Let $n \geq 1$. Then  $\wt{f_b}^n(\wt{\gamma})$ has endpoint $p$ and links with $h^i(\wt{\gamma})$ for all $1 \leq i \leq n$.}

\

\textit{Proof of Claim 1:} By our choice of $\wt{f_b}$ and $\wt{\gamma}$, we know the claim holds for $n=1$. Let $n \geq 2$.  Inductively, suppose that $\wt{f_b}^{n-1}(\wt{\gamma})$ has endpoint $p$ and links with $h^i(\wt{\gamma})$ for all $1 \leq i \leq n-1$. 
Let $I$ be the interval in $\partial \Hy$ that connects the endpoints of $\wt{\partial X_{ab}}$ and does not contain $p$, oriented from the repelling fixed point of $h$ to the attracting fixed point. We will use interval notation when speaking about connected subsets of $I$. Now $\wt{f_b}$ extends continuously to a homeomorphism of $\partial \Hy$, which we will also denote by $\wt{f_b}$, and because $\wt{f_b}$ fixes the endpoints of $\wt{\partial X_{ab}}$, this extension restricts to a homeomorphism of $I$. Let $z$ be the endpoint of $\wt{\gamma}$ in $I$, and let $x \in \partial I$ be the attracting fixed point of $h$. Because  $\wt{f_b}^{n-1}(\wt{\gamma})$   links with $h^i(\wt{\gamma})$ for all $1 \leq i \leq n-1$ and has endpoint $p$, we have
\begin{equation} \label{intervalcontainment} (\wt{f_b}^{n-1}(z), x] \subseteq (h^i(z), x] \hspace{10pt} \text{ for all } \hspace{10pt} 0 \leq i \leq n-1.
\end{equation}
Since $\wt{f_b}(\wt{\gamma})$ has endpoint $p$ and links with $h(\wt{\gamma})$, it must be that $\wt{f_b}(z) \in (hz, x]$. 
It follows from this, the fact that $\wt{f_b}$ and $h$ fix $x$, that $\wt{f_b}$ and $h$ commute by uniqueness of map lifting, and (\ref{intervalcontainment}),
that for all $0 \leq i \leq n-1$
\begin{equation} \label{intervalcontainment2} \wt{f_b}^n(z) =\wt{f_b}^{n-1}(\wt{f_b}(z))  \in \wt{f_b}^{n-1}(h (z), x]=h(\wt{f_b}^{n-1}(z), x] \subseteq h (h^i(z), x] =(h^{i+1}(z), x].
\end{equation}
Because $\wt{f_b}$ fixes $p$, we have $\wt{f_b}^n(p)=p$. This combined with (\ref{intervalcontainment2}) implies that $\wt{f_b}^n(\wt{\gamma})$ links with $h^{i+1}(\wt{\gamma})$ for all $0 \leq i \leq n-1$, proving  Claim 1.

\ 

By Claim 1, after replacing $\wt{f_b}^{n}(\wt{\gamma})$ with the geodesic connecting its endpoints, the images of $\wt{f_b}^n(\wt{\gamma})$ and $\wt{\gamma}$  in $(\overline{\mathbb{H}}^2-\{x,y\})/ \langle h \rangle$ intersect each other at least $n$ times, and all these intersections have the same sign. 
Now apply Equation (\ref{eq:annulusDistance}) to see that 
\[d_{\partial X_{ab}}(\gamma, \phi(b^n)\gamma) \geq n+1.\]
 It follows that
\begin{align} \label{1} d_{\partial X_{ab}}(\mu, \phi(b^n)\mu) &\geq d_{\partial X_{ab}}(\gamma, \phi(b^n)\gamma) -d_{\partial X_{ab}}(\mu, \gamma) -d_{\partial X_{ab}}(\phi(b^n)\mu, \phi(b^n) \gamma) \nonumber \\ 
&\geq  n+1-2d_{\partial X_{ab}}(\mu,\gamma)
\end{align}
Lemma \ref{lemann} says that $|d_{\partial X_{ab}}(\mu, \phi(a^nb^n) \mu) - d_{\partial X_{ab}}(\mu, \phi(b^n)\mu)| \leq 4$. This together with Equation (\ref{1})  implies that  \[d_{\partial X_{ab}}(\mu, \phi(a^nb^n) \mu)) \succ n. \]
From this and that fact that  $\mathcal{C}(\partial X_{ab})$ is quasi-isometric to $\mathbb{R}$ it is immediate that \linebreak  $(\pi_{\partial X_{ab}}(\phi(a^nb^n) \mu))_{n\in \mathbb{N}}$ has a subsequence converging to a point in $\partial \mathcal{C}(\partial X_{ab})$. 
\end{proof}

\begin{proposition} \label{PropDiffLim}
The sequences $(\phi(a^{n}))_{n \in \mathbb{N}}$ and $(\phi(a^{n}b^{n}))_{n \in \mathbb{N}}$ do not converge to the same point in $\Mod(S) \cup \partial \Mod(S)$.  
\end{proposition}

\begin{proof}
After passing to a subsequence, we may assume that  $(\phi(a^{n}))_{n \in \mathbb{N}}$ and $(\phi(a^{n}b^{n}))_{n \in \mathbb{N}}$ converge to points $p$ and $q$ respectively in  $\Mod(S) \cup \partial \Mod(S)$ and, by Lemma \ref{lemlingrowth}, that $(\pi_{\partial X_{ab}} (\phi(a^nb^n) \mu))_{n \in \mathbb{N}}$  converges to a point in $\partial \C(\partial X_{ab})$. 
Lemmas \ref{lemmaGeneralLin} and  \ref{lemlingrowth} imply that $q$ is in $\partial \Mod(S)$. Say $\displaystyle{q=\sum_{Y \subseteq S} c^q_Y \lambda_Y^q}$, where $c^q_Y \geq 0$ and $\lambda_Y^q \in \partial \mathcal{C}(Y)$ for all $Y \subseteq S$. Then Lemmas \ref{lemmaGeneralLin} and \ref{lemlingrowth} also imply that $c_{\partial X_{ab}}^q >0$. 

Now if $p$ were in $\Mod(S)$, then we would be done since clearly then $p \neq q$. So we will assume that $p \in  \partial \Mod(S)$, and let $\displaystyle{p=\sum_{Y \subseteq S} c^p_Y \lambda_Y^p}$.
Now observe that by Lemma \ref{lemann} and Theorem \ref{thrm:boundProjmultiCurve}
\[d_{\partial X_{ab}}(\mu, \phi(a^{n})\mu) \leq d_{\partial X_{ab}}(\mu, \mu) +4 \leq 5.\]
Thus, $(\pi_{\partial X_{ab}}(\phi(a^{n})\mu))_{n \in \mathbb{N}}$ does not limit to a point on $\partial \C(\partial X_{ab})$.
So by definition of the topology of $\Mod(S) \cup \partial \Mod(S)$, it must be that $c^{p}_{\partial X_{ab}}=0$. Since $c^{q}_{\partial X_{ab}} >0$, we see that $p \neq q$, which completes the proof. 
\end{proof}

\section{Koberda RAAGs} \label{sec:KRAAG}

In this section we complete the proof of Theorem \ref{theoremCLMKintro}. Following this, we will discuss how to use similar techniques to prove a large class of Koberda embeddings do not extend. 

\begin{figure}
 \begin{minipage}{.65\textwidth}
  \centering

 \begin{tikzpicture}[scale=3]
 \begin{scope}
\draw[fill=gray!50] (0,0) circle (1);
\clip (0,0) circle (1);

\hglinewhite{130}{175} 
\hglinewhite{5}{50} 
\hglinewhite{30}{150} 
\hglinewhite{170}{178}
\hglinewhite{2}{10}
\hglinewhite{177}{179}
\hglinewhite{3}{1}
\draw [fill=white,white] (-1,0) rectangle (1,-1); 

\hgline{90}{270} 

\hgline{30}{150} 

\hgline{130}{175} 
\hgline{5}{50} 

\hgline{170}{178}
\hgline{2}{10}

\hgline{177}{179}
\hgline{3}{1}

\hgline{160}{270}
\hgline{140}{270}
\hgline{40}{-90}
\hgline{12}{-90}
\hgline{20}{-25} 
\end{scope}

\node[fill=black, circle, scale=.4] at (0,-1){} ;
\node[fill=black, circle, scale=.4] at (-0.939693, 0.34202){}; 

\node[fill=black, circle, scale=.4] at(-0.766044, 0.642788){}; 
\node[fill=black, circle, scale=.4] at (0.978148, 0.207912){}; 
\node[fill=black, circle, scale=.4] at(0.766044, 0.642788){}; 

\node[fill=black, circle, scale=.4] at (-0.642788, 0.766044){}; 
\node[fill=black, circle, scale=.4] at (0.642788, 0.766044){}; 
\node[fill=black, circle, scale=.4] at (0.906308, -0.422618){}; 
\node[fill=black, circle, scale=.4] at (-1,0){}; 
\node[fill=black, circle, scale=.4] at (-0.866025, 0.5){}; 

\node at (0,.15) {\small $\wt{A}$};
\node at (-1.1,0) {\small $\wt{\eta}$};
\node at (0,-1.1) {\small $p$};
\node at (-.99,0.54) {\small $\wt{\alpha}$};
\node at (-.65,.90) {\small $\wt{\beta}$};
\node at (.65,.90) {\small $h\wt{\beta}$};
\node at (-1.05,.35){\small $\wt{\gamma}$};
\node at (1.05,-.5){\small $h(\wt{\gamma})$};
\node at (-1.05,.73){\small $\wt{\phi(b)}^{c_3}(\wt{\gamma})$};
\node at (1.25,.70){\small $\wt{\phi(a)}^{c_2}\wt{\phi(b)}^{c_3}(\wt{\gamma})$};
\node at (1.6,.25){\small $ \wt{\phi(b)}^{c_1}\wt{\phi(a)}^{c_2}\wt{\phi(b)}^{c_3}(\wt{\gamma})$};

\end{tikzpicture} 

\end{minipage}
\hfill
\begin{minipage}{.55\textwidth}
\centering
   \def\svgwidth{190pt}
 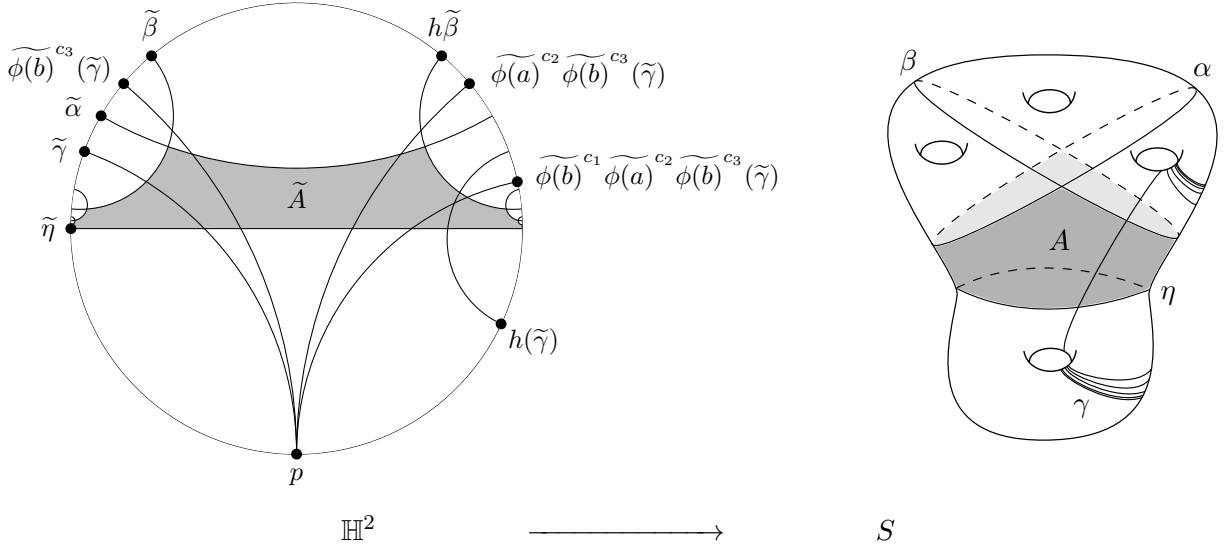
\end{minipage}%
\vspace{7pt}
\hspace{100pt} 
$\textcolor{white}{spacespacespacespacespace} \mathbb{H}^2 \hspace{54pt} \xrightarrow{\hspace*{2.4cm}} \hspace{54pt} S$
\captionsetup{width=.82\linewidth}
 \captionof{figure}{Curves $\alpha$, $\beta$, and $\eta$, bounding an annulus $A$, and simple bi-infinite geodesic $\gamma$ on surface $S$, and the universal cover $\Hy$ of $S$ as in Lemma \ref{lemmaKg}.}
  \label{fig:KUnivCover}
\end{figure}

Let  $\alpha$ and $\beta$ be the pair of intersecting curves on $S=\mathbb{H}^2/\Lambda$ depicted in Figure \ref{fig:KUnivCover}. Let $\Gamma$ be the graph with $V(\Gamma)=\{a,b\}$ and no edges. For sufficiently large $N$, Theorem \ref{Kemb} says that the homomorphism
\[ \phi: A(\Gamma) \rightarrow \Mod(S) \hspace{10pt} \text{ defined by } \phi(a)= {T_\alpha}^N \text{ and } \phi(b)= {T_\beta}^N\]
is injective, where $T_\alpha$ and $T_\beta$ denote  Dehn twists about  $\alpha$ and $\beta$ respectively.  Throughout this section, we let $\mu$ be a fixed marking on $S$.   
Equip $A(\Gamma)$ with an HHS structure. 
 
  In this section we prove the following theorem, which will complete the proof of Theorem \ref{theoremCLMKintro}. 
 
 \begin{theorem} \label{thrmmainKK}   There exists $g \in A(\Gamma)$ such that 
the sequences $(a^n)_{n \in \mathbb{N}}$ and $(a^ng^n)_{n \in \mathbb{N}}$ converge to the same point in $\partial A(\Gamma)$, but $(\phi(a^n))_{n \in \mathbb{N}}$ and $(\phi(a^ng^n))_{n \in \mathbb{N}}$ do not converge to the same point in $\partial \Mod(S)$. 
\end{theorem}
 
 As a step towards proving Theorem \ref{thrmmainKK}, we prove the following lemma in which we construct $g \in A(\Gamma)$.  

\begin{lemma} \label{lemmaKg}There exist constants $A \geq 1$ and $B \geq 0$  and a word $g \in A(\Gamma)$ such that for all $n \geq 1$ we have $d_\eta(\mu, \phi(a^ng^n)\mu) \stackrel{A,B}{\succ} n$, where $\eta$ is the  curve shown in Figure \ref{fig:KUnivCover}.  Consequently, after passing to a subsequence, $(\pi_\eta(\phi(a^ng^n)\mu))_{n \in \mathbb{N}}$ converges to a point in $\partial \mathcal{C}(\eta)$. 
\end{lemma} 

\begin{proof} We will prove that there exist constants $c_1, c_2,c_3$ such that $g=b^{c_1}a^{c_2}b^{c_3}$ has the desired properties.

Let $A$ be the annulus in Figure \ref{fig:KUnivCover}. Let $\wt{A}$ be a component of the preimage of $A$ in $\mathbb{H}^2$. Let $\wt{\beta}$ be a component of the preimage of $\beta$ such that a segment of $\wt{\beta}$ is in the boundary of $\wt{A}$, and let $\wt{\eta}$ denote the component of the preimage of $\eta$ in the boundary of $\wt{A}$. 
Let $h \in \Lambda$ be a primitive isometry with axis $\wt{\eta}$. Let $\wt{\alpha}$ be the component of the preimage of $\alpha$ that links with $\wt{\beta}$ and $h(\wt{\beta})$ and contains a segment that is in the boundary of $\wt{A}$. 

Let $Y_\alpha$ be the component of $S-\alpha$ that contain $\eta$. To simplify arguments,  we let $\phi(a)$  denote a representative in its isotopy class that fixes $Y_\alpha$ pointwise.
Let $\wt{\phi(a)}:\mathbb{H}^2\rightarrow \mathbb{H}^2$ be the lift of $\phi(a)$ that fixes some point on $\wt{\alpha}$. 
Similarly define $Y_\beta$ to be the component of $S-\beta$ containing $\eta$, choose a representative in the isotopy class of $\phi(b)$ that fixes $Y_\beta$ pointwise, and  let $\wt{\phi(b)}$ be the lift of $\phi(b)$ that fixes some point on $\wt{\beta}$.
It then follows that 
\[\wt{\phi(a)} = \mathds{1} \text{ on } \wt{Y_\alpha} \hspace{10pt} \text{ and } \hspace{10pt} \wt{\phi(b)}=\mathds{1} \text { on } \wt{Y_\beta},\] 
where for $i\in \{\alpha, \beta\}$ we let $\wt{Y_i}$ denote the component of the preimage of $Y_i$ in $\mathbb{H}^2$ whose boundary contains $\wt{i}$. Observe that for $i \in \{a, b\}$ we have that $\wt{\phi(i)}$ fixes the endpoints of $\wt{\eta}$.

Choose a geodesic $\wt{\gamma}$ in $\mathbb{H}^2$ that links with both $\wt{\beta}$ and $\wt{\eta}$ and maps to a simple bi-infinite geodesic in $S$. Further, suppose that $\wt{\gamma} \cap \wt{Y_\alpha}\cap \wt{Y_\beta}$ is an infinite ray, and let $p$ denote its endpoint on $\partial \mathbb{H}^2$. For example, take $\gamma$ to be the simple bi-infinite geodesic in $S$ with one end spiraling around a curve essential in $Y_\alpha \cap Y_\beta$ and the other end spiraling around a curve essential in $S-Y_\beta$ as in Figure \ref{fig:KUnivCover}, and  take $\wt{\gamma}$ to be an appropriate component of the preimage of $\gamma$.
Observe that $\wt{\phi(a)}$ and $\wt{\phi(b)}$ must fix $p$. 

 Now choose  $c_3 \in \mathbb{Z}$ so that $\wt{\phi(b)}^{c_3}(\wt{\gamma})$ links with $\wt{\alpha}$. Then choose $c_2 \in \mathbb{Z}$ so that $\wt{\phi(a)}^{c_2} \wt{\phi(b)}^{c_3} (\wt{\gamma})$ links with $h(\wt{\beta})$. Finally, choose $c_1 \in \mathbb{Z}$ so that $\wt{\phi(b)}^{c_1}\wt{\phi(a)}^{c_2} \wt{\phi(b)}^{c_3} (\wt{\gamma})$ links with $h(\wt{\gamma})$. See Figure \ref{fig:KUnivCover}. 

To simplify notation, define
\[g=b^{c_1}a^{c_2}b^{c_3} \in A(\Gamma) \hspace{10pt} \text { and } \hspace{10pt} \wt{\phi(g)}=\wt{\phi(b)}^{c_1}\wt{\phi(a)}^{c_2} \wt{\phi(b)}^{c_3}.\]
As in  Lemma \ref{lemlingrowth},  we have that $\wt{\phi(g)}^n(\wt{\gamma})$ has endpoint $p$ and  links with $h^i(\wt{\gamma})$ for all $1 \leq i \leq n$, implying that $d_\eta(\gamma, \phi(g^n) \gamma) \geq n+1$. 
 It follows that 
\begin{equation} \label{K1} d_{\eta}(\mu, \phi(g^n) \mu) \geq d_{\eta}(\gamma, \phi(g^n) \gamma)-d_{\eta}(\mu, \gamma) -d_{\eta}(\phi(g^n) \mu, \phi(g^n) \gamma ) \geq n+1 -2d_{\eta}(\mu, \gamma).
\end{equation}
Now Lemma \ref{lemann} says that $|d_{\eta}(\mu, \phi(a^ng^n) \mu) - d_{\eta}(\mu, \phi(g^n)\mu)| \leq 4$. This together with Equation (\ref{K1})  implies that  $d_{\eta}(\mu, \phi(a^ng^n )\mu) \succ n. $
From this and that fact that  $\mathcal{C}(\eta)$ is quasi-isometric to $\mathbb{R}$, it is immediate that  $(\pi_{\eta}(\phi(a^ng^n) \mu))_{n \in \mathbb{N}}$ has a subsequence converging to a point in $\partial \mathcal{C}(\eta)$. 
\end{proof}

We can now prove Theorem \ref{thrmmainKK}. 

\begin{proof} [Proof of Theorem \ref{thrmmainKK}] \label{proof:thrmmainKK}
Let $g\in A(\Gamma)$ be as in Lemma \ref{lemmaKg}. By the discussion in Section \ref{subsec:BoundaryHHS}, to show that $(a^n)_{n \in \mathbb{N}}$ and $(a^ng^n)_{n \in \mathbb{N}}$ converge to the same point $\partial A(\Gamma)$ it is enough to show that they  converge to the same point in $\partial_GX$, where $X$ is the Cayley graph of $A(\Gamma)$. Now the  Gromov product 
\[(a^i, a^jg^j)_1=(a^i, a^j(b^{c_1}a^{c_2}b^{c_3})^j)_1= \min(i,j) \rightarrow \infty \text{ as } i, j \rightarrow \infty. \]
Therefore $\displaystyle{\lim_{n \rightarrow \infty} a^n=\lim_{n\rightarrow \infty} a^ng^n}$ in $\partial_G X$, as desired. 

To finish this proof, we mimic the proof of Proposition \ref{PropDiffLim}. Replacing $b$ with $g$, and $\partial X_{ab}$ with $\eta$, and Lemma \ref{lemlingrowth} with Lemma \ref{lemmaKg}, we find that
$(\phi(a^{n}))_{n \in \mathbb{N}}$ and $(\phi(a^{n}g^{n}))_{n \in \mathbb{N}}$ do not converge to the same point in $\partial \Mod(S)$.
\end{proof}

Our techniques used to prove Theorem \ref{thrmmainKK} can be used to prove a more general statement on non-existence of boundary maps for right-angled Artin groups that are not necessarily free groups. To prove this more general statement, one needs to understand HHS structures for all right-angled Artin groups. 
In the following theorem, by a \textit{standard HHS structure on $A(\Gamma)$}, we mean one induced by a factor system generated by a rich family of subgraphs of $\Gamma$. 
We refer the reader to \cite{BHS}, specifically Proposition 8.3 and Remark 13.2, for details and to \cite{DHS} for a general description of the corresponding HHS boundary. In the proof of the following theorem, we freely use definitions and notations used in \cite{BHS} and \cite{DHS}.  

\begin{theorem} \label{theorem:NofillNoextend}
Let $\{\alpha_1, \ldots, \alpha_k\}$ be any collection pairwise distinct of curves in $S$. Let  $\Gamma$ be the graph with  $V(\Gamma)=\{s_1, \ldots, s_k\}$ and $s_is_j $ an edge in $\Gamma$ if and only  and $i(\alpha_i,\alpha_j)=0$.  Give $A(\Gamma)$ a standard HHS structure, or if $A(\Gamma)$ is a free group, any HHS structure. 
If there exists distinct intersecting curves $\alpha_i$ and $\alpha_j$ that do not fill $S$, then any corresponding Koberda embedding $\phi: A(\Gamma) \rightarrow \Mod(S)$ does not extend continuously to a map $\partial A(\Gamma) \rightarrow \Mod(S)$. 
\end{theorem}

 \begin{proof}
  Consider the subgraph $\Lambda$ of $\Gamma$ with $V(\Lambda)=\{s_i, s_j\}$. 
Contained in the Salvetti complex $S_\Gamma$ associated to $\Gamma$ there is a subcomplex that is the Salvetti complex associated to $A(\Lambda)$. We let $\wt{S_{\Lambda}}$ denote the lift of this subcomplex  to the universal cover $\wt{S_\Gamma}$ of $S_\Gamma$ that contains 1.  Let $\mathcal{R}$ be a rich family of induced subgraphs of $\Gamma$, and let $\mathcal{F}$ be the corresponding factor system in $\wt{S_\Gamma}$. Lemma 8.4 of \cite{BHS} tells us that 
 \[\mathcal{F}'=\{ F \cap \wt{S_\Lambda}: F \in \mathcal{F} \} \]
 is a factor system in $\wt{S_\Lambda}$. Associating $A(\Gamma)$ and $A(\Lambda)$ with $\wt{S_\Gamma}$ and  $\wt{S_\Lambda}$ respectively,  we equip each with the HHS structures corresponding to their respective factor systems. We first argue that  the  inclusion map $A(\Lambda) \rightarrow A(\Gamma)$ extends continuously to a map $\partial A(\Lambda) \rightarrow \partial A(\Gamma)$. If $\phi$  extends continuously to a map $\partial A(\Gamma)\rightarrow \partial \Mod(S)$, it will follow that $A(\Lambda) \rightarrow \Mod(S)$ extends continuously to a map $\partial A(\Lambda) \rightarrow \partial \Mod(S)$; we will show that this is impossible.

First, consider $A(\Lambda) \rightarrow A(\Gamma)$. Given $U \in \mathcal{F}'$ such that $U$ is not a 0-cube, define $\pi(U)$ to be the parallelism class of the $\subseteq$-minimal $F \in \mathcal{F}$ such that $U=F \cap \wt{S_\Lambda}$. Observe that $U$ and $V$ are nested (respectively orthogonal) if and only if $\pi(U)$ and $\pi(V)$ are nested (respectively orthogonal). This together with Lemma 10.11 of \cite{DHS} implies that $(A(\Lambda) \rightarrow A(\Gamma), \pi)$ is a hieromorphism. Theorem 5.6 of \cite{DHS} gives a condition guaranteeing that a hieromorphism extends continuously. In our case, if the following claims are true, we can apply Theorem 5.6 to conclude that $A(\Lambda) \rightarrow A(\Gamma)$ extends continuously. 

\

\noindent \textbf{Claim 1:} $\pi$ is injective. 

\noindent \textit{Proof of Claim 1:} Suppose $U, V \in \mathcal{F}'$ and $\pi(U)=\pi(V)$. Then $\pi(U) \sqsubseteq \pi(V)$ and $\pi(V) \sqsubseteq \pi(U)$. Thus, $U \subseteq V$ and $V \subseteq U$, implying $U=V$, as desired. 

\

\noindent \textbf{Claim 2:} If $[F] \in \ol{\mathcal{F}}$  is not a class of 0-cubes and  there exists no $U \in \mathcal{F}'$ satisfying $ \pi(U)=[F]$, then  $\text{diam}_{\widehat{C}F}(\pi_F(\wt{S_\Lambda}))$ is bounded above uniformly for some (any) $F \in [F]$.

\noindent \textit{Proof of Claim 2:} Let $[F]\in \overline{\mathcal{F}}$ be as in Claim 1.  First, suppose there exists $F \in [F]$ such that $F \cap \wt{S_\Lambda} \neq \emptyset$. By Lemma 8.5 in \cite{BHS}, we have  $\textgoth{g}_F(\wt{S_\Lambda}) \subseteq F \cap \wt{S_\Lambda}$. If $F \cap \wt{S_\Lambda}$ is a 0-cube, then $\text{diam}_{\widehat{C}F}(\pi_F(\wt{S_\Lambda})) \leq 1$, so the claim holds. Otherwise, there must exists $\ol{F} \in \mathcal{F}$ such that  $\ol{F} \subsetneq F$ and $\ol{F} \cap \wt{S_\Lambda}=F\cap \wt{S_\Lambda}$. It follows that $C\ol{F}$ is coned off in $\widehat{C}F$ and that $\textgoth{g}_F(\wt{S_\Lambda}) \subseteq \ol{F}$.  This implies that $\text{diam}_{\widehat{C}F}(\pi_F(\wt{S_\Lambda})) \leq 4$. 

Now assume  $F \cap \wt{S_\Lambda} =\emptyset$ for all $F \in [F]$. Choose $g \in A(\Gamma)$ and $\Gamma' \in \mathcal{R}$ so that $g\wt{S_{\Gamma'}} \in [F]$. An argument like that in the proof of Proposition 8.3 of \cite{BHS} shows that 
\begin{equation} \label{eq:gateSLambda} \textgoth{g}_{g\wt{S_{\Gamma'}}}(\wt{S_\Lambda}) =g(\wt{S}_{\Gamma' \cap \Lambda \cap \text{Lk}g}) \subseteq g(\wt{S}_{\Gamma'\cap \text{Lk}g}),  
\end{equation}
where $\text{Lk}g$ denotes the link of $g$.  Now if $\Gamma'\cap \Lambda \cap \text{Lk}g=\emptyset$, then $\textgoth{g}_{g\wt{S_{\Gamma'}}}(\wt{S_\Lambda}) =\{g\}$, implying that $\text{diam}_{\widehat{C}(g\wt{S_{\Gamma'}})}(\pi_{g\wt{S_{\Gamma'}}}(\wt{S_\Lambda}))\leq 1$. Assume then that $\Gamma'\cap \Lambda \cap\text{Lk}g \neq\emptyset$. Then by definition of $\mathcal{R}$ and $\mathcal{F}$, we have that $\Gamma'\cap \text{Lk}g \in \mathcal{R}$ and $g(\wt{S}_{\Gamma'\cap \text{Lk}g}) \in \mathcal{F}-\{0\text{-cubes}\}$. If $g(\wt{S}_{\Gamma'\cap \text{Lk}g})$ is not a proper subcomplex of $g\wt{S_{\Gamma'}}$, then $\Gamma' \subseteq \text{Lk}g$, implying  that $\wt{S_{\Gamma'}}$ is parallel to $g\wt{S_{\Gamma'}}$ (see Lemma 2.4 in \cite{BHS}). But this cannot be because $\wt{S_{\Gamma'}} \cap \wt{S_\Lambda} =\wt{S}_{\Gamma' \cap \Lambda}\neq \emptyset$ and no factor parallel to $g\wt{S_{\Gamma'}}$ intersects $\wt{S_\Lambda}$ non-trivially. Therefore, $g(\wt{S}_{\Gamma'\cap \text{Lk}g})$ must be a proper subcomplex of $g\wt{S_{\Gamma'}}$. Thus,  $Cg(\wt{S}_{\Gamma'\cap \text{Lk}g})$ is coned off in $\widehat{C}(g\wt{S_{\Gamma'}})$. This together with (\ref{eq:gateSLambda}) implies that $\text{diam}_{\widehat{C}(g\wt{S_{\Gamma'}})}(\pi_{g\wt{S_{\Gamma'}}}(\wt{S_\Lambda}))\leq 4$,  completing the proof of  Claim 2. 

\
  
We now argue that $A(\Lambda) \rightarrow \Mod(S)$ does not extend continuously to a map $\partial A(\Lambda) \rightarrow \partial \Mod(S)$. Let $\eta$ denote a geodesic representative of an essential boundary component of a small regular neighborhood of $\alpha_i \cup \alpha_j$. Using the proof techniques of Lemma \ref{lemmaKg}, we can construct  $g \in A(\Lambda)$  so that 
 $d_\eta(\mu, \phi(s_i^ng^n)\mu)$ grows linearly in $n$.  For later convenience, we construct $g$ so that when written in reduced form, the first letter of $g$  is $s_j^{\pm1}$.   As in Proposition \ref{PropDiffLim}, we see that the  sequences $(\phi(s_i^n))$ and $(\phi(s_i^ng^n))$ do not converge to the same point in $\Mod(S) \cup \partial \Mod(S)$. Now observe that $(s_i^n)$ and $(s_i^ng^n)$ converge to the same point in $\partial_G A(\Lambda)$. Therefore, by the discussion in Section \ref{subsec:BoundaryHHS}, $(s_i^n)$ and $(s_i^ng^n)$  converge to the same point in $\partial A(\Lambda)$. 
 We have now established that
 $A(\Lambda) \rightarrow \Mod(S)$ does not extend continuously to a map $\partial A(\Lambda) \rightarrow \partial \Mod(S)$. Therefore, $A(\Gamma) \rightarrow \Mod(S)$ does not extend continuously when $A(\Gamma)$ is equipped with a standard HHS structure. 
 
 Now suppose $A(\Gamma)$ is a free group equipped with any HHS structure. Then by the discussion in Section \ref{subsec:BoundaryHHS}, because $(s_i^n)$ and $(s_i^ng^n)$ converge to the same point in $\partial_G A(\Gamma)$, we have that $(s_i^n)$ and $(s_i^ng^n)$ converge to the same point in $\partial A(\Gamma)$. Because $(\phi(s_i^n))$ and $(\phi(s_i^ng^n))$ do not converge to the same point in $\partial \Mod(S)$, it follows that $A(\Gamma) \rightarrow \Mod(S)$ does not extend continuously. 
 \end{proof}

\section{Existence of boundary maps for some free groups} \label{section:existence}

In this section, we show that a class of embeddings of free groups in $\Mod(S)$, which include a class of Koberda embeddings and a class of Clay, Leininger, and Mangahas embeddings,  extend continuously.  

Throughout this section, let $\Gamma$ be the graph with $V(\Gamma)=\{s_1, \ldots, s_k\}$ and no edges, and let $A(\Gamma)$ denote the corresponding right-angled Artin group (a rank $k$ free group). Equip $A(\Gamma)$ with an HHS structure. 
Let $\{X_1, \ldots, X_k\}$ be a collection of pairwise distinct, pairwise overlapping, and pairwise filling  subsurfaces of $S$ and $\{f_1,\ldots, f_k\}$ a collection of mapping class such that $f_i$ is fully supported on $X_i$. 
Let $\mu$ be a fixed marking on $S$. 
The main theorem of this section is the following, which implies the remaining direction of Theorem \ref{theorem:IntroExistenceFill} in  the introduction. 

\begin{theorem} \label{theorem:existenceFill} Let $A(\Gamma)$ be the rank $k$ free group equipped with any HHS structure. Let $\{X_1, \ldots, X_k\}$ a collection of pairwise distinct, pairwise overlapping, and pairwise filling  subsurfaces of $S$, and $\{f_1,\ldots, f_k\}$ a collection of mapping class such that $f_i$ is fully supported on $X_i$.  There
 exists a $C>0$ such that if  $\tau_{X_i}(f_i) \geq C$ for all $i$, then  the homomorphism \[\phi: A(\Gamma) \rightarrow \Mod(S) \hspace{10pt} \text{ defined by } \hspace{10pt} \phi(s_i)= f_i  \text{ for all } i \] is a quasi-isometric embedding and extends continuously to a map $ \partial A(\Gamma) \rightarrow \partial \Mod(S)$. 
\end{theorem}

We emphasize the arguments we will use to establish that $\phi$ is a quasi-isometric embedding are essentially the same as those used by Clay, Leininger, and Mangahas to prove Theorem \ref{CLMemb}. In particular, when the the $X_i$ are all non-annular, that $\phi$ is a quasi-isometric embedding is Theorem \ref{CLMemb}.  To prove Theorem \ref{theorem:existenceFill}, we require the following proposition.

\begin{proposition} \label{prop:bigProjSylSub_PO}  There exists $K >0$ such that the following holds. For each $1 \leq i \leq k$, assume $\tau_{X_i}(f_i) \geq 2K$. Let $\phi:A(\Gamma) \rightarrow \Mod(S)$ be the homomorphism defined by $\phi(s_i)=f_i$ for all $i$.  Consider $g_1\ldots g_k \in A(\Gamma)$, where for each $i$ we have $g_i=x_i^{e_i}$ for some $x_i \in \{s_1^{\pm 1}, \ldots, s_k^{\pm 1} \}$ and $e_i >0$, and  $x_i \neq x_{i+1}$, and $x_1^{e_1} \ldots x_k^{e_k}$ is a reduced word.  Let $Y_i$ be the subsurface of $S$ that fully supports $\phi(x_i)$. Then
\begin{enumerate} \item For each  $1 \leq i \leq k$, we have  $d_{\phi(g_1\ldots g_{i-1}) Y_i}(\mu, \phi(g_1\ldots g_k) \mu) \geq Ke_i$, 
\item  For all $1 \leq i<j \leq k$, we have 
$\phi(g_1\ldots g_{i-1}) Y_i \prec \phi(g_1 \ldots g_{j-1}) Y_{j}$, where $\prec$ denotes the partial order on  $\Omega(K,\mu, \phi(g_1\ldots g_k) \mu)$, and 
\item The homomorphism $\phi:A(\Gamma) \rightarrow \Mod(S)$ is a quasi-isometric embedding. 
\end{enumerate}
\end{proposition}

\begin{proof}Define $K=K_0+20+2\max\{d_{X_i}(\mu, \partial X_j): 1 \leq i,j \leq k \text{ and } i\neq j \}$, where $K_0$ is maximum of the constants in Theorem 6.12 of \cite{MMII} and Theorem \ref{bgit}. 
Statements (1) and (2) of this proposition are essentially Theorem 5.2 in \cite{CLM}. 
The difference is that Theorem 5.2 does not allow for the homomorphism to send a generator to a power of a Dehn twist. The only obstruction to Theorem 5.2 holding for homomorphisms $\phi$ of this type is the following. Suppose $X_i$ is the subsurface that fully supports $\phi(s_i)$ and let $\sigma\in A(\Gamma)$ be a non-empty word in letters commuting with $s_i$, not including $s_i$. If $X_i$ is non-annular, then $d_{X_i}(\phi(\sigma)\mu', \mu'')=d_{X_i}(\mu', \mu'')$ for any markings $\mu', \mu''$. This not necessarily true if $X_i$ is an annulus. However, this issue does not arise for us because $A(\Gamma)$ a free group implies no such $\sigma$ exists. Thus, the arguments used to prove Theorem 5.2 in \cite{CLM} also prove our Statements (1) and (2).  The proof of our Statement (3) is the same as the proof in \cite{CLM} of Theorem \ref{CLMemb}, using our Statements (1)  instead of their Theorem 5.2. 
\end{proof}

The proof of the next lemma is essentially contained in the proof of Theorem 6.1 in \cite{CLM}. We include a proof here for completeness. 
\begin{lemma} \label{lemma:SylSubsNeighborsOnGeos}
Let $\phi: A(\Gamma) \rightarrow \Mod(S)$, $g_1\ldots g_k \in A(\Gamma)$, and $Y_i$ be as in Proposition \ref{prop:bigProjSylSub_PO}. Let $\mathcal{G}$ be a geodesic in $\C(S)$ with one end in $\pi_S(\mu)$ and one end in $\pi_S(\phi(g_1\ldots g_k) \mu)$. Then for each $1\leq i \leq k$, there exists a curve $\gamma_i$ on $\mathcal{G}$ such that $\pi_{\phi(g_1\ldots g_{i-1}) Y_i}(\gamma_i) = \emptyset$. If $|i-j| \geq 3$ and $\gamma_i $ and $\gamma_j$ are two such curves, then $\gamma_i \neq \gamma_j$. 
\end{lemma}

\begin{proof}
Fix $1 \leq i \leq k$. By way of contradiction, suppose for all curves $v$ on $\mathcal{G}$, we have $\pi_{\phi(g_1\ldots g_{i-1})Y_i}(v) \neq \emptyset$. Then  Theorem \ref{bgit} and Theorem \ref{thrm:boundProjmultiCurve} together imply that 
\begin{align*} d_{\phi(g_1\ldots g_{i-1}) Y_i}(\mu, \phi(g_1\ldots g_k) \mu) \leq 4+K_0.
\end{align*}
But Proposition \ref{prop:bigProjSylSub_PO} says $d_{\phi(g_1\ldots g_{i-1}) Y_i}(\mu, \phi(g_1\ldots g_k) \mu) \geq K >K_0+4$, a contradiction. Thus, there must exists a curve $\gamma_i$ on $\mathcal{G}$ such that $\pi_{\phi(g_1\ldots g_{i-1}) Y_i}(\gamma_i) = \emptyset$, as desired. Note that this implies that $\gamma_i$ and $\partial \phi(g_1\ldots g_{i-1}) Y_i$ form a multicurve. 

Now consider $\gamma_i$ and $\gamma_j$, where  $1 \leq i <j \leq k$ and $|i-j| \geq 3$. We will show that $\gamma_i$ and $\gamma_j$ are distinct curves. To the contrary, suppose $\gamma_i=\gamma_j$. 
Because of the filling assumption on $\{X_1, \ldots, X_k\}$, the pair of subsurfaces $ Y_{i+1}$ and $Y_{i+2}$  fill $S$. Thus,  $\phi(g_1\ldots g_{i+1})Y_{i+1}= \phi(g_1\ldots g_{i})Y_{i+1}$ and $ \phi(g_1\ldots g_{i+1}) Y_{i+2}$ are also a pair of subsurfaces that fill $S$. 
Thus, it must be that $\pi_{\phi(g_1\ldots g_{n-1})Y_n}(\gamma) \neq \emptyset$ for some $n \in \{i+1, i+2\}$. In any case, $i < n < j$. 
 
In the remainder of this proof, to simplify notation, for each $\ell$ we define $\overline{Y_\ell}= \phi(g_1 \ldots g_{\ell-1}) Y_\ell$.
 By Proposition \ref{prop:bigProjSylSub_PO}, we have 
\[ \ol{Y_i} \prec  \ol{Y_n}  \prec  \ol{Y_j} ,\]
where $\prec$ is the partial order on $\Omega(K, \mu, \phi(g_1 \ldots g_k)\mu)$. 
In particular, these three subsurfaces are pairwise overlapping. This together with the assumption that $\gamma_i =\gamma_j$ and Theorem \ref{thrm:boundProjmultiCurve} implies that 
\[d_{\ol{Y_n}}(\partial \ol{Y_i}, \partial \ol{Y_j})  \leq  d_{\ol{Y_n}}(\partial \ol{Y_i}, \gamma_i)+d_{\ol{Y_n}}(\gamma_j, \partial \ol{ Y_j}) \leq 2+2=4.
\]
It follows from this and the definition of  $\prec$ that
\[d_{\ol{Y_n}}(\mu, \phi(g_1 \ldots g_k) \mu) 
\leq d_{\ol{Y_n}}(\mu, \partial \ol{Y_i}) +d_{\ol{Y_n}}(\partial \ol{Y_i},\partial \ol{Y_j}) +d_{\ol{Y_n}}(\partial \ol{Y_j},\phi(g_1 \ldots g_k) \mu) \leq 4 +4+4 =12.\]
But this cannot be, because $d_{\ol{Y_n}}(\mu, \phi(g_1 \ldots g_k) \mu) \geq K \geq 20$ by Proposition \ref{prop:bigProjSylSub_PO}.  Therefore, $\gamma_i$ and $\gamma_j$ are distinct curves. 
\end{proof}

We have now developed the tools we will need to prove Theorem \ref{theorem:existenceFill}.
\begin{proof}[Proof of Theorem \ref{theorem:existenceFill}] 
Define $C=2K$, where $K$ is as in Proposition \ref{prop:bigProjSylSub_PO} and for each $1 \leq i \leq k$, assume that $\tau_{X_i}(f_i) \geq C$.  By Proposition \ref{prop:bigProjSylSub_PO}, $\phi$ is a quasi-isometric embedding. 

Let $X$ denote the Cayley graph of $A(\Gamma)$. Choose $x \in \partial_GX$. 
Let $\gamma$ be the infinite geodesic ray in $X$ based at 1 limiting to a point $x$ in $\partial_G X$. 
We think of $\gamma$ as an infinite word of the form $y_1y_2y_3\ldots$, where each  $y_i \in \{s_1^{\pm1}, \ldots, s_k^{\pm 1}\}$ and the word $y_1 y_2 \ldots y_i$ is a reduced word for all $i$. By construction, the sequence  $(y_1\ldots y_n)$ in converges to $x$ in $X \cup \partial_G X$. 
Let $(h_n)$ be another sequence in $A(\Gamma)$ that converges to $x$ in $X \cup \partial_GX$. We will show that $(\phi(h_n))$ and $(\phi(y_1\ldots y_n))$ converge to the same point in $\partial \Mod(S)$. By the discussion in Section \ref{subsec:BoundaryHHS}, this will prove the theorem. 
  We will consider two case: 
(1) There does not exist $N\geq 1$ such that $y_i=y_N $ for all $i \geq N$, and (2) such an $N$ exists.
In both cases, we will assume each $h_n$ is written in the form $h_n=g_{n,1}\ldots g_{n,N(n)}$, where for all $i$ we have $g_{n,i}=x_{n,i}^{e_{n,i}}$ for some $e_{n,i} >0$ and $x_{n,i} \in \{s_1^{\pm 1} , \ldots, s_k^{\pm 1} \}$ satisfying $x_{n,i} \neq x_{n, i+1}$, and $x_{n,1}^{e,1} \ldots x_{n,N(n)}^{e_{n,N(n)}}$ is a reduced word.

\textbf{Case 1:} Suppose there does not exist $N\geq 1$ such that   $y_i=y_N$ for all $i \geq N$. Then we can think of $\gamma$ as an infinite word of the form $g_1g_2g_3\ldots$, where  $g_i=x_i^{e_i}$ for some $e_i >0$ and  $x_i \in \{s_1^{\pm 1}, \ldots, s_k^{\pm 1} \}$ satisfying $x_i \neq x_{i+1}$, and $x_1^{e_1} \ldots x_i^{e_i}$ is a reduced word for all $i$. Define $Y_i$ to be the subsurface that fully supports $\phi(x_i)$. 
For short, we let $\overline{Y_i}$ denote $\phi(g_1 \ldots g_{i-1})Y_i$. 

Because $(h_n)$ and $(y_1\ldots y_n)$ converge to the same point in $\partial_G X$ and $X$ is a tree, $h_n$ and $y_1\ldots y_n$ must agree on longer and longer initial segments as $n\rightarrow \infty$. In particular, 
 given $L \geq 1$, there exists $M$ such that for all $n \geq M$, we have $g_{n,1}\ldots g_{n,L}=g_1\ldots g_L$.
Consider $n \geq M$ and $k \geq e_1+\cdots +e_L$. Choose a curve $\beta \in \text{base}(\mu)$. Given $\sigma \in A(\Gamma)$, let $\mathcal{G}(\sigma)$ denote some choice of geodesic  in $\C(S)$ with endpoints $\beta$ and $\phi(\sigma)\beta$. 
By Lemma \ref{lemma:SylSubsNeighborsOnGeos}, for all $1 \leq i \leq L$ there exist curves $\gamma_i$ and  $\gamma_i'$ on $\mathcal{G}(y_1\ldots y_k)$ and $\mathcal{G}(h_n)$ respectively such that $\pi_{\ol{Y_i}}(\gamma_i) =\emptyset$ and $\pi_{\ol{Y_i}}(\gamma_i') =\emptyset$. Observe that
\[d_S(\gamma_i, \partial \ol{Y_i})\leq 1 \hspace{10pt} \text{ and }  \hspace{10pt} d_S(\gamma_i', \partial \ol{Y_i}) \leq 1.\]

Choose $\gamma_r$ to be the curve in  $\{\gamma_i: 1 \leq i \leq L\}$  closest to $\phi(y_1\ldots y_k)\beta$. Lemma \ref{lemma:SylSubsNeighborsOnGeos} tells us that if $|i-j| \geq 3$, then $\gamma_i \neq \gamma_j$. So necessarily $d_S(\beta, \gamma_r) \geq L/3$.  Thus, the Gromov product, computed in $\C(S)$, 
\begin{align*}
(\phi(y_1\ldots y_k)\beta, \phi(h_n) \beta)_\beta 
&=\frac{1}{2}\bigg[d_S(\beta, \phi(y_1\ldots y_k)\beta)+d_S(\beta, \phi(h_n)\beta)-d_S(\phi(y_1\ldots y_k)\beta, \phi(h_n) \beta) \bigg]\\
& \geq \frac{1}{2}\bigg[ d_S(\beta, \gamma_r)+d_S(\gamma_r, \phi(y_1\ldots y_k)\beta)+d_S(\beta, \gamma_r') +d_S(\gamma_r', \phi(h_n)\beta) -\\
&\hspace{5pt}  \bigg(d_S(\phi(y_1\ldots y_k)\beta, \gamma_r)+d_S(\gamma_r, \partial \ol{Y_r}) +d_S(\partial \ol{Y_r}, \gamma_r') +d_S(\gamma_r', \phi(h_n)\beta)\bigg) \bigg] \\
&\geq \frac{1}{2}\bigg[d_S(\beta, \gamma_r)+d_S(\beta, \gamma_r') -2\bigg] \\
&\geq \frac{1}{2}(L/3-2).
\end{align*}
It follows that 
\begin{equation} \label{eq:GromovProductInfty} \displaystyle{\liminf_{k, n \rightarrow \infty} (\phi(y_1\ldots y_k)\beta, \phi(h_n)\beta))_{\beta} =\infty}.
\end{equation} 
 Because $(h_n)$ is an arbitrary sequence converging to $x$, we could have taken it to be $(y_1\ldots y_n)$. Thus, Equation (\ref{eq:GromovProductInfty}) tells us two things: (1) $(\phi(y_1 \ldots y_n)\mu)$ converges to a point in $\partial \C(S)$, and (2) $(\phi(y_1\ldots y_n)\mu)$ and $(\phi(h_n) \mu)$ converge to the same point in  $\partial \C(S)$. By definition of the topology on $\Mod(S)\cup \partial \Mod(S)$, this tells us that $(\phi(y_1\ldots y_n))$ and $(\phi(h_n))$ converge to the same point in $\partial \Mod(S)$. 

\textbf{Case 2:}  Assume there exists $N\geq 1$ such that $y_i=y_N $ for all $i \geq N$.  Corollary 6.2 in \cite{DHS} tells us that the action of $\Mod(S)$ by left multiplication extends to an action of $\Mod(S)$ on $\Mod(S)  \cup \partial \Mod(S)$ by homeomorphisms.
Consequently, if we can show that $(\phi((y_1\ldots y_{N-1})^{-1}h_n))_{n \in \mathbb{N}}$ and $(\phi(y_N\ldots y_n))_{n \in \mathbb{N}}$ converge to the same point in $\partial \Mod(S)$, then $(\phi(h_n))_{n \in \mathbb{N}}$ and $(\phi(y_1\ldots y_n))_{n \in \mathbb{N}}$ must converge to the same point in $\partial \Mod(S)$.  Furthermore, $((y_1\ldots y_{N-1})^{-1}h_n)_{n \in \mathbb{N}}$ and $(y_N\ldots y_n)_{n \in \mathbb{N}}$ converge to the same point in $\partial_G X$. Thus, without loss of generality we  assume $N=1$. 
By our assumption, $y_1\ldots y_n=y_1^n$ for all $n$. 

Let $Y$ be the subsurface that fully supports $\phi(y_1)$ and let $\partial Y=\{\beta_1, \ldots, \beta_\ell\}$. Then
\begin{equation} \label{eq:limY}
\lim_{n\rightarrow \infty} \frac{d_Y(\mu, \phi(y_1^n) \mu)}{n} >0 \hspace{10pt}  \text{ and } \hspace{10pt}\pi_Y(\phi(y_1^n)\mu) \rightarrow \lambda_Y \text{ for some } \lambda_Y \in \partial \C(Y).
\end{equation} 
Further observe that for all $i$
\begin{equation} \label{eq:limbetai}
\lim_{n\rightarrow \infty} \frac{d_{\beta_i}(\mu, \phi(y_1^n) \mu)}{n} \geq0.
\end{equation}
If (\ref{eq:limbetai}) is an equality, let $\lambda_i$ be any point in $\partial \C(\beta_i)$. Otherwise, define $\lambda_i \in \partial \C(\beta_i)$ to be $\displaystyle{\lim_{n\rightarrow \infty} \pi_{\beta_i}(\phi(y_1^n)\mu).}$ For all subsurfaces $W$ disjoint from $Y$ and not an annulus with core curve in $\partial Y$,  Lemma \ref{lemann} and Theorem \ref{thrm:boundProjmultiCurve} imply that $d_W(\mu, \phi(y_1^n)\mu) \leq d_W(\mu, \mu)+4 \leq 6 $. 
Consequently, \[\lim_{n \rightarrow \infty} \phi(y_1^n) =c_Y\lambda_Y +\sum_{i=1}^\ell c_i\lambda_i, \text{ where } c_Y+\sum_{i=1}^\ell c_i=1 \text{ and } \frac{c_i}{c_Y}=\lim_{n \rightarrow \infty} \frac{d_{\beta_i}(\mu, \phi(y_1^n)\mu)}{d_Y(\mu, \phi(y_1^n)\mu)}.\] 

 Because  $(h_n)$ and $(y_1^n)$ converge to the same point in $\partial_G X$, given any $L \geq 1$, for all sufficiently large $n$ we have $x_{n,1}=y_1$ and $e_{n,1} \geq L$. 
 So by removing finitely many initial terms from $h_n$, for convenience we may assume that $g_{n,1}=y_1^{e_{n,1}}$ for all $n$. Observe that $e_{n,1} \rightarrow \infty$ as $n\rightarrow \infty$. 
 It is immediate from this and the definition of the topology of $\Mod(S) \cup \partial \Mod(S)$ that $\displaystyle{\lim_{n \rightarrow \infty} \phi(g_{n,1}) =\lim_{n \rightarrow \infty} \phi(y_1^n)}$.  
Thus,  to finish the proof, we must show $\displaystyle{\lim_{n \rightarrow \infty} \phi(g_{n,1})=\lim_{n \rightarrow \infty} \phi(h_n)}$.  
By passing to subsequences, we may assume that either $N(n)=1$ for all $n$ or $N(n)\geq 2$ for all $n$. If the former holds, then $h_n=g_{n,1}$, and we are done. Assume then that  $N(n) \geq 2$ for all $n$.
To proceed, we require the following claims. 

\noindent \textbf{Claim 1:} $d_Y(\phi(g_{n,1})\mu, \phi(h_n)\mu)$ is bounded above, independent of $n$.

\noindent \textbf{Claim 2:} Let $W$ be a subsurface that is disjoint from $Y$. Then $d_W(\phi(g_{n,1})\mu, \phi(h_n)\mu)$ is bounded above, independent of $n$.

We postpone the proofs of these claims and for now assume they are true. First, observe that Claim 1 and (\ref{eq:limY}) imply that $\pi_Y(\phi(h_n)\mu) \rightarrow \lambda_Y.$ If Inequality (\ref{eq:limbetai}) is strict, then Claim 2 implies that $\pi_{\beta_i}(\phi(h_n)\mu) \rightarrow \lambda_i$. Further observe that Claims 1 and 2 imply that for all $W$ disjoint from $Y$
\[\lim_{n\rightarrow \infty} \frac{d_W(\mu, \phi(g_{n,1})\mu)}{d_Y(\mu, \phi(g_{n,1})\mu)} =\frac{\displaystyle{ \lim_{n\rightarrow \infty} \frac{d_W(\mu, \phi(g_{n,1})\mu)}{e_{n,1}}}}{ \displaystyle{\lim_{n\rightarrow \infty}  \frac{d_Y(\mu, \phi(g_{n,1})\mu)}{e_{n,1}}}}
=\frac{ \displaystyle{\lim_{n\rightarrow \infty} }\frac{d_W(\mu, \phi(h_n)\mu)}{e_{n,1}}}{\displaystyle{\lim_{n\rightarrow \infty}}  \frac{d_Y(\mu, \phi(h_n)\mu)}{e_{n,1}}}
=\lim_{n\rightarrow \infty} \frac{d_W(\mu, \phi(h_n)\mu)}{d_Y(\mu, \phi(h_n)\mu)}.\]
It follows that $\displaystyle{\lim_{n \rightarrow \infty} \phi(g_{n,1})=\lim_{n \rightarrow \infty} \phi(h_n)}$ as desired. 
 
To finish the proof, we will now prove Claims 1 and 2.  For each $n$, let $Z_n$ denote the subsurface that fully supports $\phi(x_{n,2})$.  

\noindent \textit{Proof of Claim 1:} Fix $n\geq 1$.
Because $Y$ fully supports $\phi(x_{n,1})$, by Proposition \ref{prop:bigProjSylSub_PO}, we know $Y \prec \phi(g_{n,1}) Z_n$, where $\prec $ denotes the partial order on $\Omega(K, \mu, \phi(h_n)\mu)$.
Thus, \linebreak $d_{Y}(\partial \phi(g_{n,1})Z_n, \phi(h_n) \mu) \leq 4$. Therefore, 
\[ d_{Y}(\phi(g_{n,1})\mu, \phi(h_n)\mu)
\leq d_{Y}(\phi(g_{n,1})\mu, \partial \phi(g_{n,1}) Z_n)+d_{Y}(\partial \phi(g_{n,1}) Z_n, \phi(h_n) \mu) 
\leq d_{Y}(\mu, \partial Z_n)+4.\]
There are finitely many possibilities for $Z_n$, so this completes the proof of Claim 1. 

\

\noindent \textit{Proof of Claim 2:} Fix $n \geq 1$. 
Because  $Y$ and $Z_n$ fill $S$ and $Y$ and $W$ are disjoint, it must be that $\pi_{Z_n}(\partial W) \neq \emptyset$. There are two cases to consider: (1) $W \pf Z_n$ and (2) $W \subsetneq Z_n$. First, suppose that $W \pf Z_n$. It then follows from Proposition \ref{prop:bigProjSylSub_PO}, Theorem \ref{thrm:boundProjmultiCurve}, and the definition of $K$ that 
\begin{align*}
 d_{Z_n}(\partial W, \phi(g_{n,2}\ldots g_{n,N(n)}) \mu) 
& \geq d_{Z_n}(\mu,\phi(g_{n,2}\ldots g_{n,N(n)})\mu)-d_{Z_n}(\partial Y, \partial W) -d_{Z_n}(\mu, \partial Y)  \\
& \geq K-2 -K/2   \geq 10. 
\end{align*}
Thus Theorem \ref{bineq} implies that $d_W(\partial Z_n, \phi(g_{n,2}\ldots g_{n,N(n)}) \mu) \leq 4$.
From this and Theorem \ref{theorem:MMLipProjection} we find that
\begin{align} \label{eq:boundedprojW}
d_W(\phi(g_{n,1})\mu, \phi(h_n)\mu) 
&= d_W(\mu, \phi(g_{n,2}\ldots g_{n,N(n)}) \mu)  \nonumber \\
& \leq d_W(\mu, \partial Z_n) +d_W(\partial Z_n, \phi(g_{n,2}\ldots g_{n,N(n)}) \mu) \nonumber \\ 
& \leq 4\max\{ d_{\wt{\mathcal{M}}(S)}(\mu, \mu_i): 1\leq i \leq k\}+4,
 \end{align}
 where $\mu_i$ is a fixed choice of marking with $\partial X_i \subseteq \text{base}(\mu_i)$  for each $1 \leq i \leq k$. This provides a uniform bound in the case that $W \pf Z_n$. 
 
 Now suppose that $W \subsetneq Z_n$. First, observe that because $Z_n$ fully supports $\phi(x_{n,2})$, the sequence $(\pi_{Z_n}(\phi(x_{n,2})^m \mu))_{m \in \mathbb{N}}$ converges to a point in $\partial \C(Z_n)$. Thus, by Corollary \ref{bgitCor} there exists a constant $M$, that depends on $W$ and $x_{n,2}$, such that $d_W(\mu, \phi(g_{n,2})\mu) \leq M$ for all $n$. Note that there are only finitely many possibilities for $x_{n,2}$, so $M$ can be chosen to be independent of $n$.  This  implies that 
 \begin{align*} d_W(\phi(g_{n,1})\mu, \phi(h_n) \mu) &\leq d_W(\mu, \phi(g_{n,2}\ldots g_{n,N(n)}) \mu)  \\
 & =d_W(\mu ,\phi(g_{n,2})\mu)+d_W(\phi(g_{n,2})\mu, \phi(g_{n,2}\ldots g_{n,N(n)})\mu) \\
 & \leq M +d_{\phi(g_{n,2})^{-1}W}(\mu, \phi(g_{n,3}\ldots g_{n, N(n)})\mu).
 \end{align*}
  Now if $N(n)=2$, then we can apply Theorem \ref{thrm:boundProjmultiCurve} to see that  \[d_{\phi(g_{n,2})^{-1}W}(\mu, \phi(g_{n,3}\ldots g_{n, N(n)})\mu)=d_{\phi(g_{n,2})^{-1}W}(\mu, \mu)\leq 2,\] and Claim 2 is established. Suppose then that $N(n) \geq 3$. Let $V_n$ denote the subsurface that fully supports $\phi(x_{n,3})$. Observe that because $\tau_{Z_n}(\phi(x_{n,2})) \geq 2K$ and $\partial Y$ and $\partial W$ form a multicurve, we have
  \begin{align*}
  d_{Z_n}(\partial \phi(g_{n,2})^{-1}W, \partial V_n) & \geq d_{Z_n}(\partial W, \partial \phi(g_{n,2})^{-1}W) -d_{Z_n}(\partial W, \partial Y) -d_{Z_n}(\mu, \partial Y)-d_{Z_n}(\mu,\partial V_n) \\
  &\geq 2K-2-K/2-K/2 > 2. 
  \end{align*}
  This together with Theorem \ref{theorem:MMLipProjection} establishes that $\partial \phi(g_{n,2})^{-1}W$ and $\partial V_n$ do not form a multicurve. Thus, $\phi(g_{n,2})^{-1} W \pf V_n$. So to bound $d_{\phi(g_{n,2})^{-1}W}(\mu, \phi(g_{n,3}\ldots g_{n, N(n)})\mu)$ from above independent of $n$, we can use the same techniques used above to bound $d_W(\mu, \phi(g_{n,2}\ldots g_{n,N(n)})\mu)$ when $W \pf Z_n$. This completes the proof of Claim 2. 
 \end{proof}

\bibliography{thesis_sources}{}
\bibliographystyle{plain}

\end{document}